\documentclass[a4paper,12pt]{amsart}

\usepackage{amsmath,amsfonts,amsthm,amsopn,cite,mathrsfs,amssymb}
\usepackage{mathtools}
\usepackage[T1]{fontenc}
\usepackage{lmodern}
\usepackage[a4paper,left=3cm,right=3cm,top=3.5cm,bottom=3.5cm]{geometry}
\usepackage{subfigure}
\usepackage{xcolor}
\usepackage{url}
\usepackage{comment}
\usepackage{tikz}
\usetikzlibrary{calc,arrows.meta,decorations.markings}

\definecolor{treeInk}{HTML}{263642}
\definecolor{tubeBlue}{HTML}{147D99}
\definecolor{tubeFill}{HTML}{DCF0F4}
\definecolor{antennaBlue}{RGB}{0,90,245}

\definecolor{sphBlue}{RGB}{15,100,221}
\definecolor{sphRed}{RGB}{204,43,37}
\definecolor{sphGreen}{RGB}{104,166,46}
\definecolor{sphGold}{RGB}{192,128,0}
\definecolor{sphFill}{RGB}{252,245,232}

\newtheorem{theorem}{Theorem}[section]
\newtheorem{proposition}[theorem]{Proposition}
\newtheorem{corollary}[theorem]{Corollary}
\newtheorem{lemma}[theorem]{Lemma}
\newtheorem{definition}[theorem]{Definition}

\theoremstyle{remark}
\newtheorem{remark}[theorem]{Remark}

\newcommand{\T}{\mathbb{T}}

\DeclareMathOperator{\osc}{osc}

\DeclareMathOperator{\inj}{inj}
\DeclareMathOperator{\supp}{supp}
\DeclareMathOperator{\spann}{span}
\newcommand{\Sd}{\mathbb{S}^d}
\newcommand{\Z}{\mathbb{Z}}
\newcommand{\R}{\mathbb{R}}
\newcommand{\M}{\mathcal{M}}

\newcommand{\dd}{\,\mathrm{d}}

\begin{document}
\begin{abstract}
A spherical \(t\)-design is a finite set of points on the \(d\)-dimensional unit sphere such that averaging any polynomial of degree at most \(t\) over these points gives its exact spherical average. We study the corresponding notion for curves, where point averages are replaced by normalized arclength integration. For fixed \(d\ge2\), the length of such curves is known to be bounded below by a positive multiple of \(t^{d-1}\) as \(t\to\infty\). We show that this asymptotic order is attained: for every \(d\ge2\), there exists a constant \(C_d>0\) such that, for every positive integer \(t\) and every prescribed length \(L\ge C_dt^{d-1}\), there is a simple spherical \(t\)-design curve of length \(L\). This establishes the analogue for curves of the theorem of Bondarenko, Radchenko, and Viazovska, which resolved the Korevaar-Meyers conjecture for spherical design points.

More generally, we prove the corresponding result on closed, connected Riemannian manifolds of dimension \(d\ge2\), where  eigenfunctions of the Laplace-Beltrami operator substitute for polynomials. The proof combines a geometric construction of curves with continuously adjustable local lengths and an application of Brouwer's degree theorem.
\end{abstract}

\title[]{Optimal asymptotic bounds for design curves on manifolds}
\author{Martin Ehler}
\address{Faculty of Mathematics \\
University of Vienna \\
Oskar-Morgenstern-Platz 1 \\
A-1090 Vienna, Austria}
\email{martin.ehler@univie.ac.at}
\maketitle

\vspace{-.5cm}

\section{Introduction}
Spherical designs provide exact integration of low-degree polynomials on the unit sphere $\Sd=\{x\in\R^{d+1} : \|x\|=1\}\subseteq\R^{d+1}$. 
The classical notion is that of spherical $t$-design points \cite{DelsarteGoethalsSeidel}, where $x_1,\dots,x_N\in\mathbb{S}^d$ satisfy 
\begin{equation*}
\frac{1}{N}\sum_{i=1}^N P(x_i)
=
\int_{\mathbb S^d} P(x)\dd\mu_{\Sd}(x)
\end{equation*}
for every polynomial $P$ of total degree at most $t$, and $\mu_{\Sd}$ is the normalized surface measure.  For fixed dimension $d$ and growing $t$, every $t$-design must have cardinality $N$ at least of order $t^d$ \cite{DelsarteGoethalsSeidel}. The existence of spherical $t$-designs was established in \cite{SeymourZaslavsky}, but no upper bound on the number of points required was provided.

Korevaar and Meyers \cite{KM} proved the existence of spherical $t$-designs $x_1,\dots,x_N\in\Sd$ with $N$ of order $t^{(d^2+d)/2}$ and conjectured the existence with $N$ of order $t^d$. This problem and actual constructions of spherical designs for distinct $t$ were the subject of extensive theoretical and numerical work by a variety of 
authors over several decades, see the websites \cite{Sloane,Womersley} and references listed there.   

Bondarenko, Radchenko, and
Viazovska \cite{BondarenkoRadchenkoViazovska} ultimately settled the
conjecture by proving
that there is a constant $C_d>0$, depending only on $d$, such that for
every positive integer $t$ and every integer $N\geq C_d\,t^d$, there exists
a spherical $t$-design on $\Sd$ consisting of exactly $N$ points. The points can be chosen well-separated \cite{BondarenkoRadchenkoViazovska2}. In the meantime, Etayo, Marzo, and Ortega-Cerd\`a \cite{Etayo1} first derived analogous results for algebraic manifolds, whereas Gariboldi and Gigante~\cite{Gariboldi1} then established optimal asymptotic bounds for design points on closed, connected, orientable Riemannian manifolds.

For $t$-design curves, the discrete sum is replaced by integration with respect to arclength. We investigate the corresponding question of optimal length for $t$-design curves \cite{EhlerGrochenig,KingMixon,Lindblad}. The passage from point configurations to the one-dimensional geometry of curves introduces substantial new geometric challenges.

\subsection{Spherical design curves}
We first consider design curves on the sphere \cite{EhlerGrochenig,Lindblad,KingMixon}.
\begin{definition}[Sphere]\label{def:designs S}
A closed rectifiable curve $\gamma:[0,L]\to  \mathbb S^d$, parametrized by arclength, is called a spherical $t$-design curve if the normalized path integral satisfies
\begin{equation*}
\frac{1}{L}\int_0^L P(\gamma(s)) d s 
=
\int_{\mathbb S^d} P(x)\dd\mu_{\Sd}(x)
\end{equation*}
for every polynomial $P$ of degree at most $t$. 
\end{definition}

Just as the cardinality measures the size of finitely many points, length is the natural measure for curves.  The basic quantitative question is how short a $t$-design curve can be. A lower bound on the order is known: every $t$-design curve must have length at least of order $t^{d-1}$ \cite{EhlerGrochenig,EhlerGrochenig2,KingMixon,Kunis}.   

This asymptotic order is attained by
constructions of $t$-design curves with length $\mathcal O(t)$ on
$\mathbb S^2$ \cite{EhlerGrochenig} and $\mathcal O(t^2)$ on
$\mathbb S^3$ \cite{Lindblad}, where the latter even yields simple
curves of every prescribed length $L\ge C_3t^2$. 

In arbitrary dimension $d \geq 2$, the inductive construction in \cite{EhlerGrochenig} yields $t$-design curves on $\Sd$ with
length $\mathcal{O}(t^{d(d-1)/2})$. As observed in \cite{Lindblad}, combining the inductive
construction with the result on $\mathbb S^3$ improves this
bound to $\mathcal{O}(t^{d(d-1)/2-1})$ for $d \geq 3$. Lindblad \cite{Lindblad2} also constructs weighted $t$-design curves of length
$\mathcal{O}(t^{d-1})$ in every dimension $d\ge2$. This weaker concept allows an additional
weight function in the arclength integral. Spherical curves of length of order $t^{d-1}$ satisfying Marcinkiewicz-Zygmund inequalities with respect 
to unweighted arclength integration are derived in \cite{EhlerGroechenigKarner}. However, these inequalities do not provide the exact integration identities required of $t$-designs.  

The existence of $t$-design curves of length $\mathcal{O}(t^{d-1})$ 
remained open for $d\ge4$. Our main contribution is to show that $t^{d-1}$ is the correct order of magnitude, and every sufficiently large prescribed length can be realized.  
\begin{theorem}
\label{thm:piecewise}
For every $d\geq2$ there is a constant $C_d>0$ such
that, for all positive integers $t$ and all real $L\geq C_dt^{d-1}$, there is a simple spherical $t$-design curve of length $L$. 
\end{theorem}
This establishes the analogue for curves of the Korevaar-Meyers conjecture.

\subsection{Design curves on closed manifolds}
Let $\mathcal M$ be a closed\footnote{compact and without boundary}, connected Riemannian manifold of dimension $d\geq2$ with normalized Riemannian measure
$\mu(\mathcal M)=1$. Denote its nonnegative Laplace-Beltrami operator by
$\Delta_{\mathcal M}$, and choose a real-valued orthonormal eigenbasis
\[
  \Delta_{\mathcal M}\varphi_j=\lambda_j^2\varphi_j,
  \qquad 0=\lambda_0<\lambda_1\leq\lambda_2\leq\cdots,
  \qquad \varphi_0\equiv 1.
\]
Following the spectral convention for manifold designs
\cite{Gariboldi1,EhlerEtayo}, set
\[
  \Pi_t(\mathcal M):=\spann_{\mathbb R}
       \{\varphi_j:\lambda_j\leq t\}.
\]
This spectral subspace is the analogue of polynomials on the sphere.

\begin{definition}\label{def:manifold-design-curve}
Let $t$ be a positive real number. A closed rectifiable curve $\gamma:[0,L]\to \M$, parametrized by
arclength, is a (spectral) \emph{$t$-design curve} on $\mathcal M$ if
\[
  \frac1L\int_0^L P(\gamma(s))\,ds
  =\int_{\mathcal M} P(x)\,\dd\mu(x),
  \qquad P\in\Pi_t(\mathcal M).
\]
\end{definition}
As noted in \cite{EhlerGrochenig}, the lower bounds for the lengths of spherical $t$-design curves established therein extend to closed, connected Riemannian manifolds. In particular, for all sufficiently large $t$, every $t$-design curve on a $d$-dimensional closed, connected Riemannian manifold $\mathcal M$ must have length at least of order $t^{d-1}$. 

Even for the torus $\T^d=\mathbb{R}^d/\mathbb{Z}^d$, questions remained open. Toroidal $t$-design curves have length at least of order $t^{d-1}$ \cite{Kunis}, and simple\footnote{A curve is called simple if it has no self-intersections} curves attaining this order have been constructed \cite{Lindblad,Kunis}. However, these curves are not contractible\footnote{A curve is called contractible if it is homotopic to a point, i.e., can be continuously deformed to a point.}, and the existence of simple, contractible toroidal $t$-design curves of length $\mathcal{O}(t^{d-1})$ remained open for $d\geq 2$.

Beyond spheres and tori, asymptotically optimal constructions are scarce. The real Grassmannian $G_{2,4}$ is treated in \cite{KarnerEhlerPedro}, and 
Lindblad addresses the projective spaces $\mathbb{RP}^2$ and
$\mathbb{RP}^3$, as well as the products of tori with $\mathbb{S}^2$, $\mathbb{S}^3$,
$\mathbb{RP}^2$, or $\mathbb{RP}^3$ \cite[Section~3]{Lindblad}. 

We extend Theorem~\ref{thm:piecewise} by showing that the order $t^{d-1}$ for length of $t$-design curves is attained on every closed, connected Riemannian
manifold $\M$ of dimension $d\geq 2$. Moreover, the curve can be chosen to
be smooth or as a geodesic cycle. Here, a \emph{geodesic chain} is a curve obtained
by joining consecutive points in a finite sequence by geodesics,
and a closed geodesic chain is called a \emph{geodesic cycle}. 
\begin{theorem}\label{thm:manifold-design-curves}
There exists $C_{\mathcal M}>0$ such that
for every real $t\geq 1$ and every $L\geq C_{\mathcal M} t^{d-1}$ there is a simple geodesic cycle that forms a contractible $t$-design curve on $\mathcal M$ of length $L$. 
\end{theorem}

The geodesic cycle can also be replaced by a smooth curve.
\begin{theorem}\label{thm:manifold-design-curves2}
There exists $C_{\mathcal M}>0$ such that
for every real $t\geq 1$ and every $L\geq C_{\mathcal M} t^{d-1}$ there is a simple, smooth, contractible $t$-design curve on $\mathcal M$ of length $L$. 
\end{theorem}

These theorems also imply Theorem \ref{thm:piecewise}, because, for every positive integer $t$,
\begin{equation*}
\Pi_t(\Sd)  \subseteq \spann_{\mathbb R}\{P|_{\Sd} : P \text{ a polynomial of total degree}\,\leq t \} \subseteq \Pi_{\sqrt{d}\, t} (\Sd),
\end{equation*}
see Section \ref{app:M as T and S}.

As we will see, constructing design curves requires several new ideas beyond those used to construct design points. The requirement that the curves be simple further complicates the construction.

\subsection{Method}
Our approach is guided by the original proof of Bondarenko, Radchenko, and Viazovska \cite{BondarenkoRadchenkoViazovska}. Broadly speaking, their method and its extension to manifolds by Etayo, Marzo, and Ortega-Cerd\`a \cite{Etayo1} and Gariboldi and Gigante \cite{Gariboldi1} rely on four ingredients:
\begin{itemize}
\item[(a)] a suitable partition of the sphere;
\item[(b)] differential equations, whose dynamics move point sets on the sphere;
\item[(c)] an estimate of local gradients or oscillation of polynomials of degree $t$; 
\item[(d)] a version of Brouwer's Degree Theorem \ref{thm:brouwer} \cite{BondarenkoRadchenkoViazovska,Brouwer}.
\end{itemize}
Choosing one point from each cell (a) of an equal-area partition
$R_1,\dots,R_N$ of $\Sd$, with $N\geq C_d t^d$, yields an
initial point set of cardinality $N$. The flow of differential equations (b) moves these points, generating a continuously parametrized family of point sets, each of cardinality $N$. Local oscillation estimates (c) establish the positivity condition required to apply Brouwer's degree theorem. The theorem then guarantees that this family contains a $t$-design. 

The partition, oscillation, and degree-theoretic ingredients
(a), (c), and (d) remain available in the setting of curves. A natural attempt to adapt part (b) to curves is to apply the flow to an initial curve, thereby obtaining a continuously parametrized family of curves. However, unlike the cardinality of point configurations, the arclength depends on the geometry of the evolving curve. We were unable to control the lengths of the curves and establish the required positivity condition simultaneously.

\begin{figure}[h]
\centering
\begin{tikzpicture}[scale=.6,
     line cap=round, line join=round,
    curve/.style={blue, line width=.7pt},
    >=stealth
]
    \foreach \k in {1,...,5}
        \draw[black!30, line width=.25pt]
            (\k,0) -- (\k,6) (0,\k) -- (6,\k);
    \draw[black, line width=.6pt] (0,0) rectangle (6,6);

    \draw[curve]
        (.5,5.5) -- (5.5,5.5)
        -- (5.5,4.5) -- (1.5,4.5)
        -- (1.5,3.5) -- (5.5,3.5)
        -- (5.5,2.5) -- (1.5,2.5)
        -- (1.5,1.5) -- (5.5,1.5)
        -- (5.5,.5) -- (.5,.5) -- cycle;

    \foreach \x in {0,...,5}
        \foreach \y in {0,...,5}
            \fill[black!50] (\x+.5,\y+.5) circle (.75pt);

    \draw[-, line width=.6pt] (2.65,0) -- (3.35,0);
    \draw[-, line width=.6pt] (2.65,6) -- (3.35,6);
    \draw[-, line width=.6pt] (0,2.65) -- (0,3.35);
    \draw[-, line width=.6pt] (6,2.65) -- (6,3.35);

\end{tikzpicture}
\hspace{.1cm}
\begin{tikzpicture}[scale=.6,
    line cap=round, line join=round,
    curve/.style={blue, line width=.7pt},
    >=stealth
]
    \foreach \k in {1,...,5}
        \draw[black!20, line width=.1pt]
            (\k,0) -- (\k,6) (0,\k) -- (6,\k);
    \draw[black, line width=.6pt] (0,0) rectangle (6,6);

    \foreach \x/\y/\angle/\turn in {
        0/5/90/1, 1/5/0/0, 2/5/0/0, 3/5/0/0, 4/5/0/0, 5/5/0/1,
        5/4/270/1, 4/4/180/0, 3/4/180/0, 2/4/180/0,
        1/4/180/-1, 1/3/270/-1,
        2/3/0/0, 3/3/0/0, 4/3/0/0, 5/3/0/1,
        5/2/270/1, 4/2/180/0, 3/2/180/0, 2/2/180/0,
        1/2/180/-1, 1/1/270/-1,
        2/1/0/0, 3/1/0/0, 4/1/0/0, 5/1/0/1,
        5/0/270/1, 4/0/180/0, 3/0/180/0, 2/0/180/0,
        1/0/180/0, 0/0/180/1,
        0/1/90/0, 0/2/90/0, 0/3/90/0, 0/4/90/0
    } {
   \pgfmathtruncatemacro{\k}{mod(\x+2*\y,3)}
\ifcase\k
    \def\weightvalue{1.125}
\or
    \def\weightvalue{1.3}
\or
    \def\weightvalue{1.7}
\fi
\pgfmathsetmacro{\bumpheight}{(\weightvalue-1)/8}

        \begin{scope}[shift={(\x+.5,\y+.5)}, rotate=\angle]
            \ifnum\turn=0
                \draw[curve]
                    (-.5,0) -- (-.25,0)
                    -- (-.25,\bumpheight) -- (-.125,\bumpheight)
                    -- (-.125,-\bumpheight) -- (0,-\bumpheight)
                    -- (0,\bumpheight) -- (.125,\bumpheight)
                    -- (.125,-\bumpheight) -- (.25,-\bumpheight)
                    -- (.25,0) -- (.5,0);
            \else
                \begin{scope}[yscale=\turn]
                    \draw[curve]
                        (-.5,0) -- (-.375,0)
                        -- (-.375,\bumpheight) -- (-.25,\bumpheight)
                        -- (-.25,0) -- (-.125,0)
                        -- (-.125,\bumpheight) -- (0,\bumpheight)
                        -- (0,0) -- (\bumpheight,0)
                        -- (\bumpheight,-.125) -- (0,-.125)
                        -- (0,-.25) -- (\bumpheight,-.25)
                        -- (\bumpheight,-.375) -- (0,-.375)
                        -- (0,-.5);
                \end{scope}
            \fi
        \end{scope}
    }

    \draw[-, line width=.6pt] (2.65,0) -- (3.35,0);
    \draw[-, line width=.6pt] (2.65,6) -- (3.35,6);
    \draw[-, line width=.6pt] (0,2.65) -- (0,3.35);
    \draw[-, line width=.6pt] (6,2.65) -- (6,3.35);

\end{tikzpicture}
\hspace{.1cm}
\begin{tikzpicture}[scale=.6,
    line cap=round, line join=round,
    curve/.style={blue, line width=.7pt},
    >=stealth
]
    \foreach \k in {1,...,5}
        \draw[black!20, line width=.1pt]
            (\k,0) -- (\k,6) (0,\k) -- (6,\k);
    \draw[black, line width=.6pt] (0,0) rectangle (6,6);

    \foreach \x/\y/\angle/\turn in {
        0/5/90/1, 1/5/0/0, 2/5/0/0, 3/5/0/0, 4/5/0/0, 5/5/0/1,
        5/4/270/1, 4/4/180/0, 3/4/180/0, 2/4/180/0,
        1/4/180/-1, 1/3/270/-1,
        2/3/0/0, 3/3/0/0, 4/3/0/0, 5/3/0/1,
        5/2/270/1, 4/2/180/0, 3/2/180/0, 2/2/180/0,
        1/2/180/-1, 1/1/270/-1,
        2/1/0/0, 3/1/0/0, 4/1/0/0, 5/1/0/1,
        5/0/270/1, 4/0/180/0, 3/0/180/0, 2/0/180/0,
        1/0/180/0, 0/0/180/1,
        0/1/90/0, 0/2/90/0, 0/3/90/0, 0/4/90/0
    } {
        \pgfmathtruncatemacro{\k}{mod(\x+2*\y,3)}
        \ifcase\k
            \def\weightvalue{1.25}
        \or
            \def\weightvalue{2}
        \or
            \def\weightvalue{3}
        \fi
        \pgfmathsetmacro{\bumpheight}{(\weightvalue-1)/8}

        \begin{scope}[shift={(\x+.5,\y+.5)}, rotate=\angle]
            \ifnum\turn=0
                \draw[curve]
                    (-.5,0) -- (-.25,0)
                    -- (-.25,\bumpheight) -- (-.125,\bumpheight)
                    -- (-.125,-\bumpheight) -- (0,-\bumpheight)
                    -- (0,\bumpheight) -- (.125,\bumpheight)
                    -- (.125,-\bumpheight) -- (.25,-\bumpheight)
                    -- (.25,0) -- (.5,0);
            \else
                \begin{scope}[yscale=\turn]
                    \draw[curve]
                        (-.5,0) -- (-.375,0)
                        -- (-.375,\bumpheight) -- (-.25,\bumpheight)
                        -- (-.25,0) -- (-.125,0)
                        -- (-.125,\bumpheight) -- (0,\bumpheight)
                        -- (0,0) -- (\bumpheight,0)
                        -- (\bumpheight,-.125) -- (0,-.125)
                        -- (0,-.25) -- (\bumpheight,-.25)
                        -- (\bumpheight,-.375) -- (0,-.375)
                        -- (0,-.5);
                \end{scope}
            \fi
        \end{scope}
    }

    \draw[-, line width=.6pt] (2.65,0) -- (3.35,0);
    \draw[-, line width=.6pt] (2.65,6) -- (3.35,6);
    \draw[-, line width=.6pt] (0,2.65) -- (0,3.35);
    \draw[-, line width=.6pt] (6,2.65) -- (6,3.35);

\end{tikzpicture}
\caption{Each of the three squares represents the torus $\T^2$ and is partitioned into $6\times 6$ smaller squares. 
(left) A simple curve runs through the partition without any detours. (center, right) Local detours are added to the initial curve that depend continuously on prescribed local arclengths.}
\label{fig:torus-global-detours}
\end{figure}
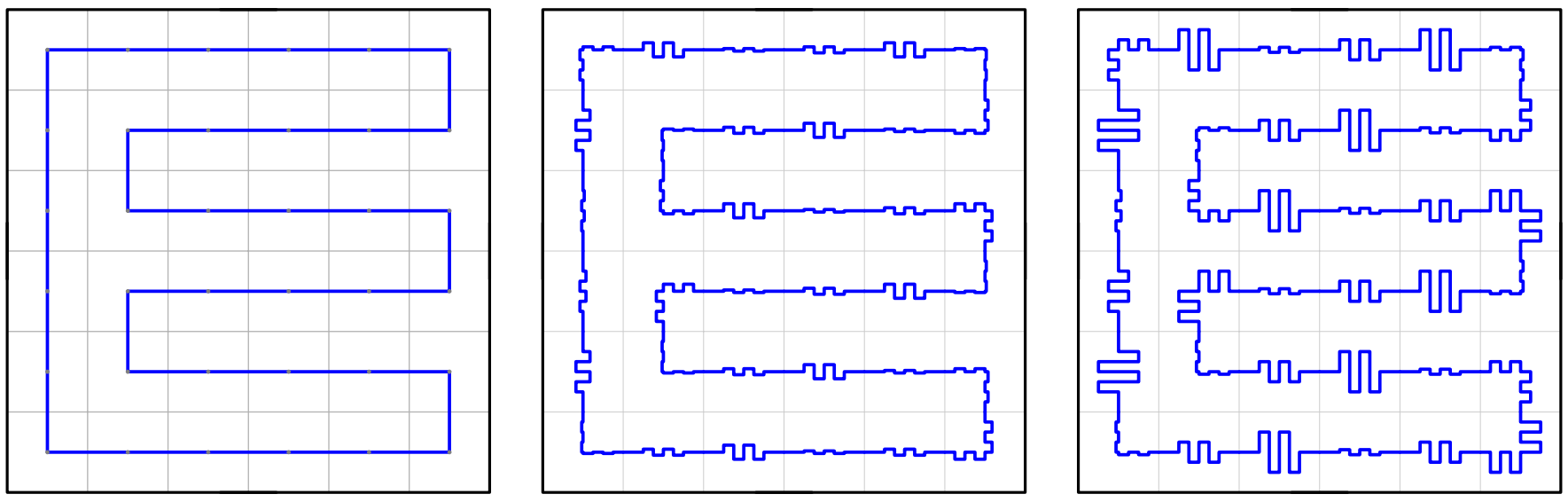

Our principal new ingredient, replacing part (b), is a geometric construction of local detours that depend continuously on prescribed local arclengths, see Figure~\ref{fig:torus-global-detours} for an illustration on the torus \(\T^2\). 
Unlike in \cite{Gariboldi1}, we do not require the manifold to be orientable. For each real $L\geq C_\M t^{d-1}$, we select an equal-area partition $R_1,\dots,R_N$ of $\M$ with $N$ of order $L^{d/(d-1)}$, and construct a simple initial curve whose length is of order $L$. The positivity condition is a parametric quadrature problem. We define parameter-dependent local quadrature weights for which the oscillation estimates ensure that the corresponding quadrature expression is positive. We then realize these weights geometrically through local detours with prescribed arclengths, obtaining a continuously parametrized family of curves with lengths of order $L$. Brouwer's degree theorem then guarantees that some member of this family integrates $\Pi_t(\M)$ exactly, hence, is a $t$-design curve. Since $\Pi_t(\M)$ contains the constant functions, the length of the selected curve is exactly $L$, rather than merely of order $L$.

\medskip 
Although not central to the paper, two observations are worth noting:

\noindent (1) \emph{Marcinkiewicz-Zygmund inequalities \cite{FilbirHielscherJahnUllrich,FilbirMhaskar,Grochenig2020,MarzoMZ,MhaskarNarcowichWard}}: At the optimal length scale $L$ of order $t^{d-1}$, the continuously parametrized family of curves from above has length $\mathcal{O}(t^{d-1})$ and satisfies Marcinkiewicz-Zygmund inequalities for all $1\leq p\leq\infty$, with uniform constants that are independent of $p$. This extends \cite{EhlerGroechenigKarner} from nonsimple spherical curves to simple, contractible curves on general closed, connected Riemannian manifolds.

\noindent (2) \emph{Logvinenko-Sereda sets \cite{MarzoCerda,Ortega,Wang}}: We formulate the interplay between local modifications and
Brouwer's degree theorem in terms of measures, with local detours
corresponding to the addition of local mass. This formulation
allows us to apply the method to geometric objects other than
curves. As another concrete manifestation of
Marcinkiewicz-Zygmund inequalities, we obtain unions of small
balls that form Logvinenko-Sereda sets and additionally satisfy
exact integration identities.

\subsection{Outline}
In Section \ref{sec:B and local}, we explain how the local detours replace part (b) in the scheme (a)--(d). The new program is executed in Section \ref{sec:polygonal-existence} on the torus, where the more elementary geometry illustrates what needs to be done on the manifold. In Section \ref{sec:manifold 1}, we derive the curves on manifolds with locally adjustable length. The oscillation estimates on the manifold are established in Section \ref{sec:osc manifold}, which verifies Theorem \ref{thm:manifold-design-curves}. A smoothing argument leads to Theorem \ref{thm:manifold-design-curves2}. Section \ref{app:M as T and S} shows that the manifold theorems imply Theorem \ref{thm:piecewise} for the sphere. The Marcinkiewicz-Zygmund inequalities and Logvinenko-Sereda sets are discussed in Section \ref{sec:further}.

\section{Local mass adjustment and Brouwer's degree theorem}\label{sec:B and local}
To describe how the mechanism of locally adjustable detours can be combined with Brouwer's degree theorem, we work in an abstract setting where curves are represented by measures and detours correspond to adding local mass. 

The construction principle is formulated in a compact metric measure space $\M$, endowed with a Borel probability measure $\mu$ that has full support. As a placeholder for $\Pi_t(\M)$, let $\Pi(\mathcal M)\subseteq\mathcal C(\mathcal M)$ be a finite-dimensional real linear subspace. Unlike in \cite{BondarenkoRadchenkoViazovska,Etayo1,Gariboldi1}, we assume that $\Pi(\mathcal M)$ contains the constant functions. Their inclusion ensures the correct total mass of the resulting measure, and hence the correct arclength when translated back to curves.

\subsection{The abstract local adjustment principle}
For every positive integer $N$, suppose we find measurable subsets $R_1,\dots,R_N\subseteq \mathcal M$ such that 
\begin{equation*}
\mathcal M=\bigcup_{i=1}^N R_i,\qquad \mu(R_i)=1/N,\qquad  \mu(R_i\cap R_j)=0, \quad i\neq j.
\end{equation*}
The sets are called an equal-area partition of $\M$. The local oscillation of $P\in\Pi(\M)$ on a compact subset $V\subseteq\M$ is defined by 
\begin{equation*}
\osc_{V}(P):=\max_{x\in V}P(x) - \min_{x\in V}P(x).
\end{equation*}
The following theorem provides the structure of our construction.
\begin{theorem}\label{thm:foundation}
Assume $N$ is a positive integer and $L,h>0$ satisfy $\frac{L}{h N}\in[\tfrac{3}{2},\tfrac{5}{2}]$. Let $R_1,\ldots,R_N$ be an equal-area partition of $\M$ and let $V_1,\dots,V_N\subseteq\M$ be compact subsets such that $R_i\subseteq V_i$. Suppose the following items hold:
\begin{itemize}
\item[(i)] There is a family of Borel measures $\{\nu_w\}_{w\in[1,3]^N}$ on $\M$, decomposed into $\nu_w=\sum_{i=1}^N\nu_{i,w}$, so that $\nu_{i,w}$ are also Borel measures satisfying
\begin{equation}\label{eq:local length}
\supp \nu_{i,w}\subseteq V_i,\qquad \nu_{i,w}(\M) = h w_i,\qquad i=1,\dots,N,
\end{equation}
and the map 
\begin{equation*}
w\mapsto
  \int_{\M}P(x)\dd \nu_w(x)
  =\sum_{i=1}^N\int_{\M}P(x)\dd \nu_{i,w}(x)
  \end{equation*}
  is continuous for all $P\in \Pi(\mathcal M)$.

\item[(ii)] The following local oscillation bound holds:
\begin{equation}\label{eq:osc finale}
\frac{12}{N}\sum_{i=1}^N \osc_{V_i}(P) < 1,\qquad \forall\, P\in\Pi(\M),\quad \|P\|_{L^1(\M)}=1.
\end{equation}
\end{itemize}
Then there is $\hat{w}\in[1,3]^N$ such that $\nu_{\hat{w}}(\M)=L$ and 
\begin{equation*}
\frac{1}{L}\int_{\M} P(x)\dd \nu_{\hat{w}}(x) = \int_{\M}P(x)\dd \mu(x),\qquad \forall\,P\in\Pi(\M).
\end{equation*}
\end{theorem}

To verify this theorem, we work with the unit $L^1(\M)$-ball of $\Pi(\M)$, 
\begin{equation*}
\Omega:=\{P\in\Pi(\mathcal M):\|P\|_{L^1(\mathcal M)}\leq1\}.
\end{equation*}
This differs from the approach in
\cite{BondarenkoRadchenkoViazovska,Etayo1,Gariboldi1},
where the norm of the gradient of $P$ is used instead.

The decisive idea is to introduce a reference measure  $\sigma_P$ by 
\begin{equation}\label{eq:measure sigmaP}
\sigma_P:=(L+hN\phi_P)\mu,\qquad \phi_P := \frac{P}{2(1+|P|)},
\end{equation}
and to replace the weight vector with a continuous function $w:\Omega\to [1,3]^N$, where each component of $w$ is 
\begin{equation}\label{eq:wi}
w_i(P):= \frac{1}{h}\sigma_P(R_i).
\end{equation}
Since $|\phi_P(x)|\leq 1/2$, the assumption $\frac{L}{h N}\in[\tfrac{3}{2},\tfrac{5}{2}]$ ensures that $\sigma_P$ is a positive measure and $w_i(P)\in[1,3]$. The mass requirement in \eqref{eq:local length} becomes
\begin{equation*}
\nu_{i,w(P)}(\M) = \sigma_P(R_i).
\end{equation*}
 The proof of Theorem \ref{thm:foundation} relies on further comparisons between $\nu_{i,w(P)}$ and $1_{R_i}\,\sigma_P$.

Let $K:\M\times \M\to\R$ denote the reproducing kernel of $\Pi(\M)$ with respect to the $L^2(\M)$-inner product, so that the reproducing property is 
\begin{equation*}
P(x)=\langle P,K(x,\cdot)\rangle_{L^2(\M)},\qquad \forall\, P\in\Pi(\M).
\end{equation*}
We define the main function 
\begin{equation}\label{eq:F}
F:\Omega\to \Pi(\M),\qquad P\mapsto \int_{\M} K(\cdot,x)\dd \nu_{w(P)}(x) - L.
\end{equation}
The term $-L$ does not appear in
\cite{BondarenkoRadchenkoViazovska,Etayo1,Gariboldi1}.
In our setting, this term, together with the inclusion of constant
functions in $\Pi(\M)$, will ensure that the total mass of $\nu_{\hat{w}}$, and hence
the length of the curve, is exactly $L$. 

By writing $F(P)$ in a fixed orthonormal basis of $\Pi(\mathcal M)$, the continuity of $w(P)$ and the integrals
in (i) imply that $F$ is continuous. 

Fubini and the reproducing property reveal, for $Q\in\Pi(\mathcal M)$ and $P\in\Omega$,
\begin{equation}\label{eq:satisfies this here}
\langle Q,F(P)\rangle_{L^2(\mathcal M)} = \int_{\M}Q(x)\dd\nu_{w(P)}(x)-L\int_{\mathcal M} Q(x)\,\dd\mu(x).
\end{equation}
If there is $P_0\in\Omega$ such that $F(P_0)=0$, then $\langle Q,F(P_0)\rangle_{L^2(\M,\mu)} = 0$, for all $Q\in \Pi(\M)$, hence, 
\begin{equation*}
\frac{1}{L}   \int_{\M}Q(x)\dd \nu_{w(P_0)}(x)  -   \int_{\M} Q(x) \dd\mu(x)  = 0,\quad \forall \,Q\in\Pi(\M),
\end{equation*}
and $\frac{1}{L}\nu_{w(P_0)}$ exactly integrates $\Pi(\M)$, so that $\hat{w}=w(P_0)$ is the sought weight vector in Theorem \ref{thm:foundation}. The choice $Q\equiv 1$ yields $\nu_{w(P_0)}(\M)=L$.

Thus, the existence problem of a measure with mass $L$ exactly integrating $\Pi(\M)$ reduces to verifying that $F$ has a zero $P_0$.

\subsection{Proof of the mechanism}
The proof of Theorem \ref{thm:foundation} is based on the following consequence of Brouwer's degree theory \cite[Theorems 1.2.6 and 1.2.9]{Brouwer}. For simplicity, we refer to it as Brouwer's degree theorem. 
\begin{theorem}[Brouwer's degree theorem]\label{thm:brouwer}
Let $H$ be a finite-dimensional real Hilbert space, let
$\|\cdot\|_*$ be any norm on $H$, and let
$\Omega=\{P\in H:\|P\|_*\leq1\}$. Suppose that
$F:\Omega\to H$ is continuous and satisfies at the boundary $\partial\Omega=\{P\in H:\|P\|_* = 1\}$,
\[
 \langle P, F(P)\rangle>0\qquad \forall \,P\in \partial\Omega.
\]
Then $F$ vanishes at a point of the interior of $\Omega$.
\end{theorem}
This has been used in \cite{BondarenkoRadchenkoViazovska,EhlerEtayo,Etayo1,Gariboldi1} for design points, and is a direct consequence of Brouwer's fixed point theorem, which may be more familiar to some readers.
\begin{proof}[Proof of Theorem \ref{thm:brouwer}]
Assume that $F$ does not have any zeros in the interior of $\Omega$. Then $f:P\mapsto -F(P)/\|F(P)\|_{*}$ is well-defined and maps the convex and compact subset $\Omega$ into itself. Brouwer's fixed point theorem implies that there is $P_0\in \Omega$ such that $f(P_0)=P_0$. This yields
\begin{equation}\label{eq:Brouwer}
0< \langle f(P_0) , f(P_0)\rangle =  \langle P_0 , f(P_0)\rangle = -\frac{1}{\|F(P_0)\|_*} \langle P_0,F(P_0)\rangle.
\end{equation}
Since $\|P_0\|_*=\|f(P_0)\|_*=1$, we have a contradiction to the positivity assumption.
\end{proof}

We verify Theorem \ref{thm:foundation}.
\begin{proof}[Proof of Theorem \ref{thm:foundation}]
For $H=\Pi(\M)$ with $L^2(\M)$-inner product and $\|\cdot\|_*$ being the $L^1(\M)$-norm, we can ensure that $F$ as defined in \eqref{eq:F} and therefore satisfying \eqref{eq:satisfies this here} has a zero $P_0$ by verifying the positivity condition
\begin{equation}\label{eq:positivity}
0 < \langle P, F(P)\rangle_{L^2(\M)} = \int_{\M}P(x)\dd \nu_{w(P)}(x)-L\int_{\mathcal M} P(x)\dd\mu(x),
\end{equation}
for $P\in\Pi(\M)$ with $\|P\|_{L^1(\M)}=1$. 

The idea is to compare $\nu_{w(P)}$ with another measure that satisfies \eqref{eq:positivity} in place of $\nu_{w(P)}$ by a sufficiently large margin. For $\phi_P=\frac{P}{2(1+|P|)}$ in \eqref{eq:measure sigmaP} with $\frac{L}{h N}\in[\tfrac{3}{2},\tfrac{5}{2}]$, the reference measure $\sigma_P:=(L+hN\phi_P)\mu$ is positive and satisfies \eqref{eq:positivity} with the margin
\begin{align*}
 \int_{\M}P(x)\dd\sigma_P(x)-L\int_{\mathcal M} P(x)\dd\mu(x)&=hN\int_{\mathcal M} P(x)\phi_P(x)\dd\mu(x)\\
  &\geq hN/4,
\end{align*}
where the inequality 
$\int_{\mathcal M} P(x)\phi_P(x)\dd\mu(x) \geq 1/4$  
 follows from Cauchy-Schwarz and $\|P\|_{L^1(\M)}=1$ via
 \begin{align*}
 1= \left\| \frac{|P|}{\sqrt{1+|P|}} \sqrt{1+|P|} \right\|^2_{L^1(\M)}
 &\leq \int_\M \frac{|P(x)|^2}{1+|P(x)|}\dd\mu(x) \int_\M (1+|P(x)|)\dd\mu(x)\\
&=4\int_\M P(x)\phi_P(x)\dd\mu(x).
 \end{align*}
 We compare $\nu_{w(P)}$ with $\sigma_P$ by
 \begin{align*}
E(P):=\int_{\M}P(x)\dd \nu_{w(P)}(x)-\int_{\mathcal M} P(x)\dd\sigma_P(x)
\end{align*}
and deduce
 \begin{equation}\label{eq:proceed from here}
 \langle P, F(P)\rangle_{L^2(\M)} =E(P)+\int_{\M}P(x)\dd \sigma_P(x)-L\int_{\mathcal M} P(x)\dd\mu(x) \geq hN/4 - |E(P)|.
 \end{equation}

According to $\nu_{w(P)} = \sum_{i=1}^N \nu_{i,w(P)}$, we localize $E(P)$ by 
\begin{equation*}
E_i(P):=\int_{\M}P(x)\dd \nu_{i,w(P)}(x)-\int_{R_i} P(x)\dd\sigma_P(x).
\end{equation*}
 The support and mass requirements \eqref{eq:local length} for $\nu_{i,w}$ with $w_i(P)$ defined in \eqref{eq:wi} lead to 
\begin{equation*}
\nu_{i,w(P)}(V_i)=\nu_{i,w(P)}(\M)= hw_i(P) = \sigma_P(R_i).
\end{equation*}
Since $R_i\subseteq V_i$, we can therefore estimate
\begin{align*}
|E_i(P)| &\leq hw_i(P) \big(\max_{x\in V_i}P(x) - \min_{x\in V_i}P(x)\big) \\
&= hw_i(P) \osc_{V_i}(P)\\
&\leq 3 h \osc_{V_i}(P),
\end{align*}
where the last inequality is due to $w_i(P)\leq 3$. Summing over $i$ and \eqref{eq:proceed from here} lead to
\begin{align*}
\langle P, F(P)\rangle_{L^2(\M)} &\geq \frac{1}{4}hN-|E(P)| \\
& \geq \frac{1}{4} hN-\frac{3h N}{N}\sum_{i=1}^N \osc_{V_i}(P) \\
&= \frac{1}{4} hN\Big(1-\frac{12}{N}\sum_{i=1}^N \osc_{V_i}(P) \Big).
\end{align*}
This term is positive due to the assumption \eqref{eq:osc finale} in Part (ii). 
\end{proof}

\begin{remark}
The strategy of first constructing a family
$\{\nu_w\}_{w\in[1,3]^N}$ provides a direct way to impose
additional properties on the resulting exact integration measure $\nu_{\hat{w}}$. Any property shared by all members of the family is inherited
by one selected from it. In the next sections, $\M$ is a closed, connected Riemannian manifold, and we construct each $\nu_w$ as the induced measure of the arclength integral of a simple curve.  
\end{remark}

\section{The torus: a model construction}
\label{sec:polygonal-existence}

We show the complete mechanism with elementary geometry and explicit constants on the torus $\T^d=\mathbb R^d/\mathbb Z^d$ for $d\geq 2$. A closed rectifiable curve on $\T^d$ of length $L>0$ is represented by a path $\gamma:[0,L]\to\mathbb{R}^d$, parametrized by
arclength, satisfying
\[
\gamma(L)-\gamma(0)\in\mathbb{Z}^d.
\]
The corresponding closed curve on $\T^d$ is obtained by
reduction modulo $\mathbb{Z}^d$. The following definition is due to \cite{Lindblad,Kunis}.
\begin{definition}[Torus]\label{def:torus}
A closed rectifiable curve $\gamma:[0,L]\to  \mathbb T^d$, parametrized by arclength, is called a toroidal $t$-design curve if the normalized path integral satisfies
\begin{equation*}
\frac{1}{L}\int_0^L e^{2\pi i \langle m,\gamma(s)\rangle} ds 
=
\int_{\mathbb T^d} e^{2\pi i \langle m,x\rangle}\dd x
\end{equation*}
for all $m\in\Z^d$ such that $\|m\|_\infty\leq t$. 
\end{definition}

Since Brouwer's degree theorem in Section \ref{sec:B and local} deals exclusively with real vector spaces, we avoid complex exponentials and instead work with the real vector space
\[
  \mathcal T_t(\mathbb T^d)
  :=\operatorname{span}_{\mathbb R}
  \left\{
    \cos(2\pi\langle m,\cdot\rangle),
    \sin(2\pi\langle m,\cdot\rangle):
    m\in\mathbb Z^d,\ \|m\|_\infty\leq t
  \right\}.
\]
Exact integration for this space is equivalent to the condition in
Definition \ref{def:torus}: taking real and imaginary parts of the
exponential identities gives exactness for the cosine and
sine functions, and conversely, combining their identities
gives exactness for the complex exponentials.

The nonnegative Laplace-Beltrami operator on $\T^d$ is $-\sum_{i=1}^d \partial_i^2$, eigenfunctions are $e^{2\pi i \langle m,x\rangle}$ with eigenvalues 
$4\pi^2 \|m\|_2^2$, for $m\in\Z^d$. The real-valued eigenfunctions are $\cos(2\pi\langle m,\cdot\rangle)$ and $\sin(2\pi\langle m,\cdot\rangle)$. Since $\|m\|_\infty\leq \|m\|_2\leq \sqrt{d}\|m\|_\infty$, we obtain
\begin{equation}\label{eq:Pi and T}
\Pi_{2\pi t} (\T^d)\subseteq \mathcal T_t(\mathbb T^d) \subseteq \Pi_{2\pi \sqrt{d} t} (\T^d),
\end{equation}
for every positive integer $t$. Consequently, every toroidal $t$-design curve is a spectral $2\pi t$-design curve in the sense of Definition \ref{def:manifold-design-curve}, while every spectral $2\pi\sqrt d\,t$-design curve is a toroidal $t$-design curve. These differences in scaling do not affect the asymptotic order in $t$.

For $1/h$ an even positive integer,

\begin{proposition}\label{lem:polygonal-local-length}
Let $d\geq 2$ and suppose that $1/h$ is an even positive integer. 
Decompose $[-\frac{1}{2},\frac{1}{2}]^d$ into $N=1/h^{d}$ axis-parallel, closed cubes $R_1,\dots,R_N$ of side-lengths $h$. Then there exists a family of simple geodesic cycles
\begin{equation*}
\{\gamma_w\}_{w\in[1,3]^N}\subseteq (-\tfrac{1}{2},\tfrac{1}{2})^d,
\end{equation*} 
parametrized by arclength, and partitioned into families of geodesic chains $\gamma_{i,w}$, for $i=1,\dots,N$, such that 
\begin{equation*}
\gamma_{i,w}\subseteq R_i,
\qquad \ell(\gamma_{i,w})=  h w_i,
\end{equation*}
and the mapping $w\mapsto 
\int_{\gamma_w}f=\sum_{i=1}^N\int_{\gamma_{i,w}}f$ is continuous, for all $f\in \mathcal C((-\tfrac{1}{2},\tfrac{1}{2})^d)$.
\end{proposition}
Since $\gamma_w\subseteq (-\tfrac{1}{2},\tfrac{1}{2})^d$, the curves are contractible in $\T^d$.
\begin{proof}
Denote the center of each cube $R_i$ by $x_i$. 
The parity of $1/h$ ensures that the nearest-neighbour graph of $x_1,\dots,x_N$ in $(-\frac{1}{2},\frac{1}{2})^d$ has a
Hamiltonian cycle, see Figure \ref{fig:Schlangenpfad} for dimension $d=2$. The case $d>2$ can be verified by induction. 
 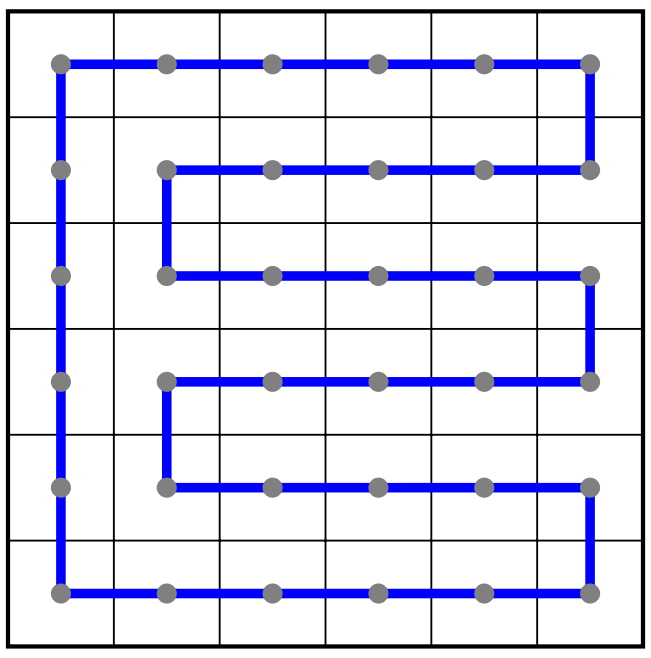
\begin{figure}
\centering
\begin{tikzpicture}[scale=.7]

  \draw[step=1,black,line width=0.35pt]
    (0,0) grid (6,6);
  \draw[black,line width=0.8pt]
    (0,0) rectangle (6,6);

  \draw[
    blue,
    line width=1.8pt,
    line cap=round,
    line join=round
  ]
    (0.5,5.5) -- (5.5,5.5)
              -- (5.5,4.5) -- (1.5,4.5)
              -- (1.5,3.5) -- (5.5,3.5)
              -- (5.5,2.5) -- (1.5,2.5)
              -- (1.5,1.5) -- (5.5,1.5)
              -- (5.5,0.5) -- (0.5,0.5)
              -- cycle;

  \foreach \x in {0.5,1.5,...,5.5}{
    \foreach \y in {0.5,1.5,...,5.5}{
      \fill[gray] (\x,\y) circle[radius=2.7pt];
    }
  }

\end{tikzpicture}
 \caption{Cuboidal partition $R_1,\dots,R_N$ of the unit cube with centers $x_i\in R_i$ that are traversed by a Hamiltonian cycle of the neighbour graph.}\label{fig:Schlangenpfad}
\end{figure}

We now have a simple, closed, piecewise linear curve. To satisfy the local length conditions within each cube $R_i$, 
say the curve enters and leaves $R_i$ in the points $p$ and $q$. We construct local detours that change continuously with $w$. 

If $p$, $x_i$, and $q$ lie on a straight line, then construct detours as in Figure \ref{fig:bumps}. 
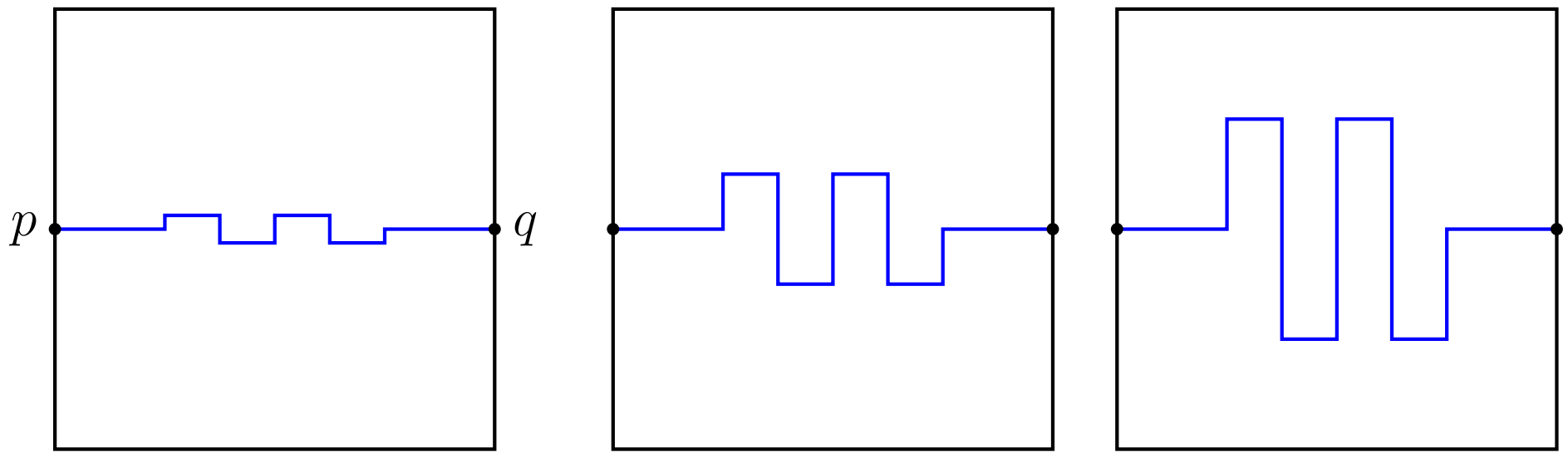
\begin{figure}
\subfigure[$w_i=5/4$]{
\begin{tikzpicture}[scale=3.5]

  \draw[thick] (0,0) rectangle (1,1);

  \draw[blue, thick]
    (0,1/2)
    -- (1/4,1/2)
    -- (1/4,17/32)
    -- (3/8,17/32)
    -- (3/8,15/32)
    -- (1/2,15/32)
    -- (1/2,17/32)
    -- (5/8,17/32)
    -- (5/8,15/32)
    -- (3/4,15/32)
    -- (3/4,1/2)
    -- (1,1/2);
  \fill (0,1/2) circle (0.4pt);
  \fill (1,1/2) circle (0.4pt);

  \node[left]  at (0,1/2) {$p$};
  \node[right] at (1,1/2) {$q$};
\end{tikzpicture}
}
\subfigure[$w_i=2$]{
\begin{tikzpicture}[scale=3.5]

  \draw[thick] (0,0) rectangle (1,1);

  \draw[blue, thick]
    (0,0.5)
    -- (0.25,0.5)
    -- (0.25,0.625)
    -- (0.375,0.625)
    -- (0.375,0.375)
    -- (0.5,0.375)
    -- (0.5,0.625)
    -- (0.625,0.625)
    -- (0.625,0.375)
    -- (0.75,0.375)
    -- (0.75,0.5)
    -- (1,0.5);
  \fill (0,1/2) circle (0.4pt);
  \fill (1,1/2) circle (0.4pt);
\end{tikzpicture}
}
\subfigure[$w_i=3$]{
\begin{tikzpicture}[scale=3.5]

  \draw[thick] (0,0) rectangle (1,1);

  \draw[blue, thick]
    (0,0.5)
    -- (0.25,0.5)
    -- (0.25,0.75)
    -- (0.375,0.75)
    -- (0.375,0.25)
    -- (0.5,0.25)
    -- (0.5,0.75)
    -- (0.625,0.75)
    -- (0.625,0.25)
    -- (0.75,0.25)
    -- (0.75,0.5)
    -- (1,0.5);
  \fill (0,1/2) circle (0.4pt);
  \fill (1,1/2) circle (0.4pt);
\end{tikzpicture}
}
\caption{In a cube $R_i$ of sidelength $h$, the bumps start and end at distance $h/4$ from $p$ and $q$. The 
height of the bumps is $\pm h(w_i-1)/8$.}\label{fig:bumps}
\end{figure}
The height of the bumps is $h(w_i-1)/8$, so that $w_i\in[1,3]$ ensures that they stay within the interior of $R_i$. If there is a corner at $x_i$, we construct detours as in Figure \ref{fig:bumps corner}. 
\begin{figure}
\subfigure[$w_i=5/4$]{
\begin{tikzpicture}[scale=3.5]

  \def\a{1/32}

  \draw[thick] (0,0) rectangle (1,1);

  \draw[blue, thick]
    (0,1/2)
    -- (1/8,1/2)
    -- (1/8,{1/2+\a})
    -- (1/4,{1/2+\a})
    -- (1/4,1/2)
    -- (3/8,1/2)
    -- (3/8,{1/2+\a})
    -- (1/2,{1/2+\a})
    -- (1/2,1/2);

  \draw[blue, thick]
    (1/2,0)
    -- (1/2,1/8)
    -- ({1/2+\a},1/8)
    -- ({1/2+\a},1/4)
    -- (1/2,1/4)
    -- (1/2,3/8)
    -- ({1/2+\a},3/8)
    -- ({1/2+\a},1/2)
    -- (1/2,1/2);

  \fill (0,1/2) circle (0.5pt);
  \fill (1/2,0) circle (0.5pt);

  \node[left]  at (0,1/2) {$p$};
  \node[below] at (1/2,0) {$q$};

\end{tikzpicture}
}
\subfigure[$w_i=2$]{
\begin{tikzpicture}[scale=3.5]

  \def\a{1/8}

  \draw[thick] (0,0) rectangle (1,1);

  \draw[blue, thick]
    (0,1/2)
    -- (1/8,1/2)
    -- (1/8,{1/2+\a})
    -- (1/4,{1/2+\a})
    -- (1/4,1/2)
    -- (3/8,1/2)
    -- (3/8,{1/2+\a})
    -- (1/2,{1/2+\a})
    -- (1/2,1/2);

  \draw[blue, thick]
    (1/2,0)
    -- (1/2,1/8)
    -- ({1/2+\a},1/8)
    -- ({1/2+\a},1/4)
    -- (1/2,1/4)
    -- (1/2,3/8)
    -- ({1/2+\a},3/8)
    -- ({1/2+\a},1/2)
    -- (1/2,1/2);

  \fill (0,1/2) circle (0.5pt);
  \fill (1/2,0) circle (0.5pt);

  \node[below] at (1/2,0) {\phantom{$q$}};

\end{tikzpicture}}
\subfigure[$w_i=3$]{
\begin{tikzpicture}[scale=3.5]

  \def\a{1/4}

  \draw[thick] (0,0) rectangle (1,1);

  \draw[blue, thick]
    (0,1/2)
    -- (1/8,1/2)
    -- (1/8,{1/2+\a})
    -- (1/4,{1/2+\a})
    -- (1/4,1/2)
    -- (3/8,1/2)
    -- (3/8,{1/2+\a})
    -- (1/2,{1/2+\a})
    -- (1/2,1/2);

  \draw[blue, thick]
    (1/2,0)
    -- (1/2,1/8)
    -- ({1/2+\a},1/8)
    -- ({1/2+\a},1/4)
    -- (1/2,1/4)
    -- (1/2,3/8)
    -- ({1/2+\a},3/8)
    -- ({1/2+\a},1/2)
    -- (1/2,1/2);

  \fill (0,1/2) circle (0.5pt);
  \fill (1/2,0) circle (0.5pt);

  \node[below] at (1/2,0) {\phantom{$q$}};

\end{tikzpicture}}
\caption{The cube $R_i$ has sidelength $h$. The bumps start and end at distance $h/8$ from $p$ and $q$. The 
height of the bumps is $h(w_i-1)/8$.}\label{fig:bumps corner}
\end{figure}
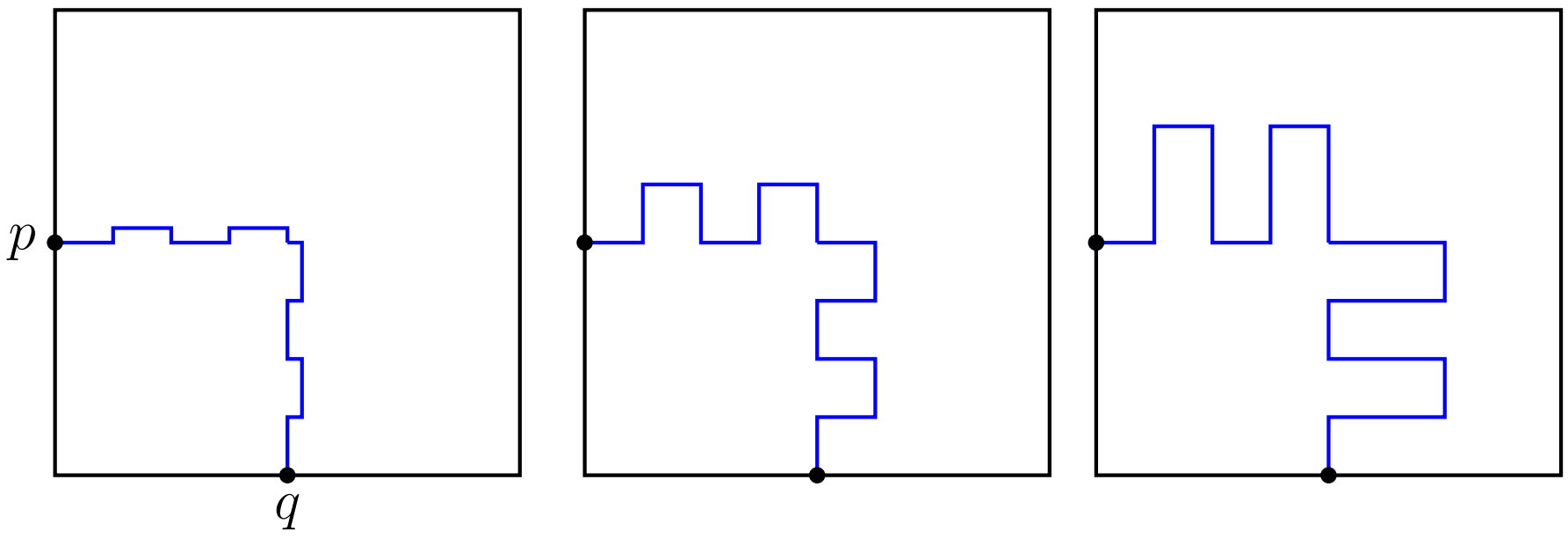
In both cases, we count the length within $R_i$ by $h$ plus $8$ times the height of the bumps. This yields
\begin{equation*}
\ell(\gamma_{i,w}) = h+ h(w_i-1) = h w_i. 
\end{equation*}
The cube $R_i$ is still entered and left at the same points $p$ and $q$, so that we obtain the simple curve with the correct local length, where assignment is simply by containment in $R_i$. Since it is piecewise linear, it can be parametrized by arclength.

Each curve $\gamma_w\subseteq (-\frac{1}{2},\frac{1}{2})^d$, see also Figure \ref{fig:torus-global-detours}, is contractible to the origin in $\T^d$ by the linear homotopy $H(\tau,s)=(1-\tau)\gamma_{w}(s)$. 

The polygon has a fixed finite number of segments with endpoints
depending continuously on $w$. 
Since segment integrals of continuous functions depend continuously
on their endpoints, $w\mapsto\int_{\gamma_w}f$ is continuous.\end{proof}

Local oscillation estimates for trigonometric polynomials on the torus are classical \cite{Zygmund}. To obtain explicit constants, we include a proof of the following estimate in Appendix \ref{app:torus}. 
\begin{proposition}[Average local oscillation]
\label{lem:torus-cubical-oscillation}
Assume $ht\leq\frac{1}{4\pi d}$, and let $P\in\mathcal T_t(\T^d)$. 
Then the averaged local oscillations satisfy
\begin{equation}\label{eq:torus-cubical-oscillation}
  \frac1N\sum_{i=1}^N\operatorname{osc}_{R_i}(P)
  \leq 4\pi d\,ht\,\|P\|_{L^1(\T^d)}.
\end{equation}
\end{proposition}

To avoid interrupting the flow of the argument, we defer the proof of the following parameter selection to Appendix~\ref{app:torus}.

\begin{lemma}[Choice of parameters]\label{lemma:parameters need balancing}
For integers $d\geq 2$ and $t\geq 1$, choose $C_d=2(64 \pi d)^{d-1}$ and $L\geq C_d t^{d-1}$ as well as $1/h= 2\lceil \frac{(L/2)^{1/(d-1)}}{2} \rceil $ with $N=1/h^d$. Then we have
\begin{equation*}
 ht \leq \frac{1}{64\pi d},\qquad \frac{L}{h N }\in[\tfrac{3}{2},\tfrac{5}{2}].
\end{equation*}
\end{lemma}

Choose $C_d$ as in Lemma \ref{lemma:parameters need balancing} and for $t\geq1$ and $L\geq C_{d}\,t^{d-1}$ choose also $N$ and $h$ as in Lemma \ref{lemma:parameters need balancing}, so that $\frac{L}{hN}\in [\frac{3}{2},\frac{5}{2}]$ and $ht \leq\frac 1{64\pi d}$ hold. 

Proposition \ref{lem:polygonal-local-length} provides a family of simple, contractible geodesic cycles $\{\gamma_w\}_{w\in[1,3]^N}$, each $\gamma_w$ partitioned into $\gamma_{i,w}$  satisfying
\[
  \gamma_{i,w}\subseteq R_i,
  \qquad \ell(\gamma_{i,w})=hw_i,
\]
as well as the continuity conditions. 
We have $ \gamma_{i,w}\subseteq R_i$, so that we may choose $V_i=R_i$. The oscillation estimate in Proposition \ref{lem:torus-cubical-oscillation}, for $P\in \mathcal{T}_t(\T^d)$ with $\|P\|_{L^1(\T^d)}=1$, yields
\begin{equation*}
\frac{12}{N}\sum_{i=1}^N\osc_{V_i}(P) \leq 48\,\pi d\, h t  \leq 48/64<1.
\end{equation*}
Thus, all assumptions of the structural Theorem \ref{thm:foundation} are satisfied, and there is $\hat{w}\in[1,3]^N$ such that 
$\gamma_{\hat{w}}$ is a toroidal $t$-design curve of length $L$. This proves Theorem \ref{thm:manifold-design-curves} in the special case of the torus.

\section{Curves on manifolds with adjustable local length}\label{sec:manifold 1}
Let $\M$ be a closed, connected Riemannian manifold of dimension $d\geq 2$. This section is dedicated to the construction of simple, contractible geodesic cycles with locally adjustable length as needed in Theorem \ref{thm:foundation}. It starts with 
\begin{itemize}
\item[(1)] a partition that induces a tree in the manifold;
\item[(2)] the tree is then approximated by a simple closed curve;
\item[(3)] this curve is approximated by a simple geodesic cycle;
\item[(4)] the geodesic cycle is modified by locally adding detours to match the local length requirements.
\end{itemize}
 Each step is carried out in a way that preserves contractibility.

We denote the Riemannian distance by $d_{\M}$ and write $B_r(x)$ for the
closed geodesic ball of radius $r$ in $\mathcal M$. Since $\M$ is closed, the injectivity radius is positive, $\inj(\M)>0$, hence, for all $x,y\in\M$, the estimate $d_\M(x,y)<\inj(\M)$ implies that $x$ and $y$ are joined by a unique minimizing geodesic segment \cite[Lemma 6.16, Corollary 6.21]{Lee}. 

We borrow the equal-measure partition theorem of Gigante and Leopardi
\cite{GiganteLeopardi}, in the manifold formulation of
\cite[Theorem 3]{Gariboldi1}.
\begin{theorem}[\cite{Gariboldi1}]\label{thm:partition}
There are $a,b>0$ such that,
for every integer $N\geq1$, there are measurable sets 
$R_1,\ldots,R_N\subseteq \M$ covering $\M = \bigcup_{i=1}^N R_i$ and points $x_i\in R_i$ such that 
\begin{align*}
  \mu(R_i)&=\frac1N,\qquad \mu(R_i\cap R_j)=0,\qquad i\neq j,\\
B_{aN^{-1/d}}(x_i) & \subseteq R_i
          \subseteq B_{bN^{-1/d}}(x_i),\qquad i=1,\dots,N.
\end{align*}
The constants $a$ and $b$ depend on the manifold $\M$, but are independent of $N$.
\end{theorem}
For $h=8bN^{-1/d}$, this partition induces a covering of $\M$ by balls $B_{h/8}(x_i)$.

\subsection{An embedded minimum spanning tree}\label{sec:embedded}
A tree with edges and vertices contained in $\M$ is called embedded if distinct edges have disjoint interiors and meet only at common vertices. 

\begin{figure}
\begin{tikzpicture}[scale=1,
  line cap=round,
  line join=round,
  >={Straight Barb[length=1.7mm,width=1.4mm]},
  tree edge/.style={black,line width=0.85pt},
  antenna/.style={antennaBlue,line width=1.25pt},
  dimension/.style={<->,black,line width=0.55pt}
]
  \coordinate (x1) at ( 0.0, 3.4);
  \coordinate (x2) at (-2.0, 0.4);
  \coordinate (x3) at ( 2.5, 0.4);
  \coordinate (x4) at (-3.5,-2.5);
  \coordinate (x5) at (-0.5,-2.5);

  \coordinate (s1) at ($(x1)+( 0.00, 0.85)$);
  \coordinate (s2) at ($(x2)+(-0.60, 0.60)$);
  \coordinate (s3) at ($(x3)+( 0.60, 0.60)$);
  \coordinate (s4) at ($(x4)+(-0.60,-0.60)$);
  \coordinate (s5) at ($(x5)+( 0.60,-0.60)$);

  \draw[tree edge]
    (x1) -- (x2)
    (x1) -- (x3)
    (x2) -- (x4)
    (x2) -- (x5);

  \node[] at (0.7,1.65) {$\Gamma$};

  \foreach \i in {1,...,5} {
    \draw[antenna] (x\i) -- (s\i);
    \fill[antennaBlue]
      (s\i) circle[radius=1.8pt];
    \filldraw[fill=black,draw=white,line width=0.45pt]
      (x\i) circle[radius=2.05pt];
  }

  \node[left=4pt]  at (x1) {$x_1$};
  \node[left=4pt]  at (x2) {$x_2$};
  \node[right=4pt] at (x3) {$x_3$};
  \node[right=4pt] at (x4) {$x_4$};
  \node[left=4pt]  at (x5) {$x_5$};

  \node[antennaBlue,above=4pt]
    at (s1) {$\sigma_{1,w_1}$};
  \node[antennaBlue,above left=3pt]
    at (s2) {$\sigma_{2,w_2}$};
  \node[antennaBlue,above right=3pt]
    at (s3) {$\sigma_{3,w_3}$};
  \node[antennaBlue,below left=3pt]
    at (s4) {$\sigma_{4,w_4}$};
  \node[antennaBlue,below right=3pt]
    at (s5) {$\sigma_{5,w_5}$};

  \draw[dimension,antennaBlue]
    ($(x1)+(0.21,0.04)$) -- ($(s1)+(0.21,-0.04)$)
    node[midway,right=3pt,font=\normalsize] {$r_1(w_1)$};

  \draw[dimension,antennaBlue]
    ($(x2)+(0.12,0.19)$) -- ($(s2)+(0.19,0.12)$);
  \node[font=\normalsize,anchor=south,antennaBlue]
    at (-1.95,1.10) {$r_2(w_2)$};

  \draw[dimension,antennaBlue]
    ($(x3)+(-0.12,0.19)$) -- ($(s3)+(-0.19,0.12)$);
  \node[font=\normalsize,anchor=south,antennaBlue]
    at (2.50,1.10) {$r_3(w_3)$};

  \draw[dimension,antennaBlue]
    ($(x4)+(0.12,-0.19)$) -- ($(s4)+(0.19,-0.12)$);
  \node[font=\normalsize,anchor=west,antennaBlue]
    at (-3.5,-3.1) {$r_4(w_4)$};

  \draw[dimension,antennaBlue]
    ($(x5)+(-0.12,-0.19)$) -- ($(s5)+(-0.19,-0.12)$);
  \node[font=\normalsize,anchor=east,antennaBlue]
    at (-0.5,-3.1) {$r_5(w_5)$};
\end{tikzpicture}
\caption{The curve $\Gamma$ traverses the tree $T$ in the order $x_1\to x_2\to x_4\to x_2\to x_5\to x_2\to x_1\to x_3\to x_1$, and the detours $\sigma_{i,w_i}$ are traversed out and back. However, this construction misses the simplicity condition. 
}\label{fig:antenna}
\end{figure}
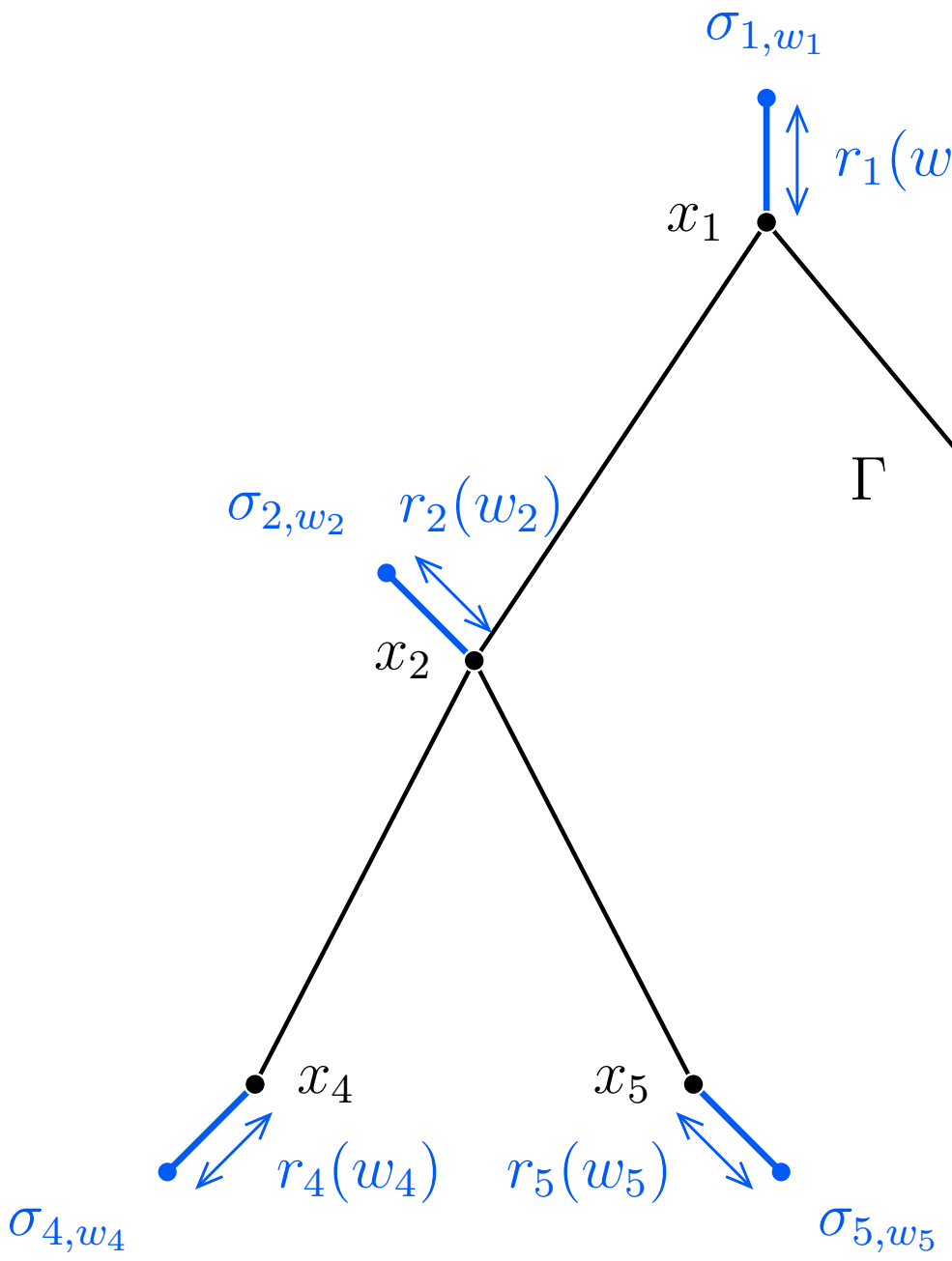

The following lemma is a Riemannian version of the classical noncrossing property of Euclidean minimum spanning trees, see \cite[Exercise 5.4.8(b)]{MatousekNesetril2008} for the planar case. 
\begin{lemma}[Embedded minimum spanning tree]\label{lemma:embedded manifold}
For $0<h<\inj(\M)$, let $x_1,\dots,x_N\in\M$ be pairwise distinct and let $G$ be the graph on these points, with an edge $e_{i,j}$ of weight $d_{\M}(x_i,x_j)$ whenever $d_{\M}(x_i,x_j)\le h/4$. If 
\begin{equation*}
\M = \bigcup_{i=1}^N B_{h/8}(x_i),
\end{equation*}
 then every minimum length spanning tree $T$ of $G$ is embedded in $\M$ by minimizing geodesics.
\end{lemma}
\begin{proof}
The covering condition ensures that $G$ is connected \cite[Lemma 4.1]{EhlerGrochenig}. Let $\widetilde{G}$ be the full graph with edge weights $d_{\M}(x_i,x_j)$. If a spanning tree of \(\widetilde G\) contains an edge longer than \(h/4\), deleting that edge produces two components. Since \(G\) is connected, some edge of \(G\) joins these components. Replacing the deleted edge by this edge strictly decreases the total length. Thus, $G$ and $\widetilde{G}$ have the same minimum length spanning trees. 

Let $T$ be a minimum length spanning tree of $G$,
and hence also of $\widetilde G$. Since $h<\inj(\M)$, every edge's geodesic realization is unique. The standard uncrossing argument, see \cite[Exercise~5.4.8(b) and its hint]{MatousekNesetril2008} for the planar case, leads to $T$ being embedded. We provide the details in Appendix \ref{app:embedded}. 
\end{proof}
The embedded tree in Lemma \ref{lemma:embedded manifold} provides the starting point for constructing curves with adjustable local length. If simplicity were not required, we could keep the tree as is and add local detours in an elementary way. Let $\Gamma$ be the piecewise geodesic cycle that traverses each edge of the tree once in each direction, and at each vertex $x_i$, we attach a geodesic segment $\sigma_{i,w_i}$ of suitably chosen length $r_i(w_i)$, traversed out and back to meet the local length requirements, see Figure \ref{fig:antenna}. Obviously, the result is not a simple curve. Constructing simple geodesic cycles with adjustable local length is more involved.

\subsection{Simple initial curves}

We invoke a general result from \cite[Theorem 2.13]{Hirsch} to approximate the tree traversal in Lemma \ref{lemma:embedded manifold} by a simple curve provided that the dimension of the manifold is at least $3$. The following is a special case of \cite[Theorem 2.13]{Hirsch}. We use $\T:=\T^1$.
\begin{lemma}[\cite{Hirsch}, approximation by embedded curves]
\label{lemma:hirsch}
Let $\mathcal M$ have dimension $d\geq3$, and let
$\gamma:\T\to\mathcal M$ be a $\mathcal C^1$-map.
For every $\varepsilon>0$, there exists a $\mathcal C^1$-embedding
$\gamma_\varepsilon:\T\to\mathcal M$ such that
\[
  d_{\mathcal M}\big(\gamma(s),\gamma_\varepsilon(s)\big)
  \leq\varepsilon,
  \qquad
   \big|
    \|\dot\gamma(s)\|-\|\dot\gamma_\varepsilon(s)\|
  \big|
  \leq\varepsilon,\qquad \forall s\in \T.
\]
In particular, $\gamma_\varepsilon$ is simple and has
positive speed everywhere.
\end{lemma}

\begin{proof}
By \cite[Chapter~2, Theorem~2.13]{Hirsch}, embeddings are
dense in $\mathcal C^1(\T,\mathcal M)$. A sufficiently close
$\mathcal C^1$ approximation satisfies both estimates.
Positive speed follows from being an embedding.
\end{proof}

In order to obtain simple geodesic chains, we may invoke an interpolation result for Jordan curves in \cite[Theorem 2.2]{Geng}. 
\begin{lemma}[Interpolation by a simple geodesic cycle]\label{lemma:interpolation}
Let $\gamma:[0,1]\to\M$ be a simple, continuous closed curve and $0<\tau_1<\dots < \tau_k<1$ points in $[0,1]$. Then for every $\varepsilon>0$ there are 
\begin{equation*}
0=s_0<s_1<\dots <s_n=1,
\end{equation*}
such that 
\begin{itemize}
\item[(i)] $s_{j+1}-s_j \leq \varepsilon$ for $j=0,\dots,n-1$;
\item[(ii)] $\{\tau_1,\dots,\tau_k\}\subseteq \{s_0,\dots,s_n\}$;
\item[(iii)]
For each $j=0,\ldots,n-1$, the points $\gamma(s_j)$ and
$\gamma(s_{j+1})$ are joined by a unique minimizing geodesic.
Traversing this geodesic at constant speed on $[s_j,s_{j+1}]$
defines an interpolating geodesic cycle $\gamma_\varepsilon$
that is simple and satisfies
$\gamma_\varepsilon(s_j)=\gamma(s_j)$;
\item[(iv)] the local length does not increase, 
\begin{equation*}
\ell(\gamma_\varepsilon|_{[s_j,s_{j+1}]}) \leq \ell(\gamma|_{[s_j,s_{j+1}]});
\end{equation*} 
\item[(v)] $d_\M(\gamma(s),\gamma_{\varepsilon}(s))\leq \varepsilon$, for all $s\in [0,1]$.
\end{itemize}
\end{lemma}
\begin{proof}
Parts (i)-(iii) are from \cite[Theorem 2.2]{Geng}. The local length estimate is obvious because we use unique minimizing geodesics. For Part (v), we observe that $\gamma$ is uniformly continuous because $[0,1]$ is compact. 
Define 
\begin{equation*}
c_\gamma(\delta):=\sup_{\substack{r,s\in[0,1]\\ |r-s|\le\delta}}
d_\M\bigl(\gamma(r),\gamma(s)\bigr).
\end{equation*} 
Uniform continuity of $\gamma$ yields $\lim_{\delta\to 0}c_\gamma(\delta)= 0$, so that we can select $0<\delta_\varepsilon\leq\varepsilon$ such that $c_\gamma(\delta_\varepsilon)<\frac{1}{2}\varepsilon$, and apply the interpolation theorem with respect to \(\delta_\varepsilon\), retaining the prescribed points. Consequently, for $s\in[s_j,s_{j+1}]$, 
\begin{align*}
d_\M\bigl(\gamma(s),\gamma_\varepsilon(s)\bigr)
&\le
d_\M\bigl(\gamma(s),\gamma(s_j)\bigr)
+
d_\M\bigl(\gamma(s_j),\gamma_\varepsilon(s)\bigr)\\
&\le
c_\gamma(\delta_\varepsilon)
+
d_\M\bigl(\gamma(s_j),\gamma(s_{j+1})\bigr)\\
&\le 2c_\gamma(\delta_\varepsilon)\leq \varepsilon.\qedhere
\end{align*}
\end{proof}

\begin{proposition}\label{prop:simple initial curve}
Suppose that $d\geq 2$. Let $R_1,\dots,R_N$ be a partition with points $x_1,\ldots,x_N$ as in Theorem \ref{thm:partition} and suppose that $h=8bN^{-1/d}<\frac{2}{3}\inj(\M)$. Then there is a simple, contractible geodesic cycle $\gamma$, whose segment traversals can be partitioned into $\gamma_1,\dots,\gamma_N$ such that 
\begin{equation*}
\gamma_i\subseteq B_{h/2}(x_i),\qquad \ell(\gamma_i)\leq h. 
\end{equation*}
\end{proposition}

\begin{figure}
\definecolor{tubeBlue}{HTML}{147D99}
\definecolor{tubeFill}{HTML}{DCF0F4}

\begin{tikzpicture}[
  scale=.9,
  line cap=round,
  line join=round,
  tree edge/.style={black,line width=0.85pt}
]
  \coordinate (x1) at ( 0, 0);
  \coordinate (x2) at (-2, 2);
  \coordinate (x3) at ( 2, 2);
  \coordinate (x4) at (-2,-2);
  \coordinate (x5) at ( 2,-2);

  \def\tuberadius{0.26}
  \pgfmathsetmacro{\tubeoffset}{\tuberadius/sqrt(2)}
  \pgfmathsetmacro{\tubecorner}{sqrt(2)*\tuberadius}

  \path[draw=tubeBlue,fill=tubeFill,line width=0.8pt]
    ({2+\tubeoffset},{2-\tubeoffset})
    arc[start angle=-45,end angle=135,radius=\tuberadius]
    [rounded corners=2pt]
    -- (0,\tubecorner)
    -- ({-2+\tubeoffset},{2+\tubeoffset})
    [sharp corners]
    arc[start angle=45,end angle=225,radius=\tuberadius]
    [rounded corners=2pt]
    -- (-\tubecorner,0)
    -- ({-2-\tubeoffset},{-2+\tubeoffset})
    [sharp corners]
    arc[start angle=135,end angle=315,radius=\tuberadius]
    [rounded corners=2pt]
    -- (0,-\tubecorner)
    -- ({2-\tubeoffset},{-2-\tubeoffset})
    [sharp corners]
    arc[start angle=225,end angle=405,radius=\tuberadius]
    [rounded corners=2pt]
-- (\tubecorner,0)
[sharp corners]
-- cycle;
  \foreach \px/\py/\ang in {
     1/ 1/45,
    -1/ 1/135,
    -1/-1/225,
     1/-1/315
  }{
    \begin{scope}[shift={(\px,\py)},rotate=\ang]
      \draw[tubeBlue,line width=0.8pt]
        (0.02,{-\tuberadius-0.09})
        -- (0.24,{-\tuberadius})
        -- (0.02,{-\tuberadius+0.09});

      \draw[tubeBlue,line width=0.8pt]
        (-0.02,{\tuberadius-0.09})
        -- (-0.24,{\tuberadius})
        -- (-0.02,{\tuberadius+0.09});
    \end{scope}
  }
  \foreach \i in {2,...,5}
    \draw[tree edge] (x1) -- (x\i);

  \foreach \i in {1,...,5}
    \fill[black] (x\i) circle[radius=2pt];

  \node[above=15pt,inner sep=0pt] at (x1) {$x_1$};
  \node[above left=5pt]        at (x2) {$x_2$};
  \node[above right=5pt]       at (x3) {$x_3$};
  \node[below left=5pt]        at (x4) {$x_4$};
  \node[below right=5pt]       at (x5) {$x_5$};
\end{tikzpicture}
\caption{The tree traversal must respect the oriented boundary of the neighborhood. For example, the traversal $x_1\to x_2\to x_1\to x_5\to x_1\to x_3\to x_1\to x_4\to x_1$ is incompatible, whereas $x_1\to x_2\to x_1\to x_4\to x_1\to x_5\to x_1\to x_3\to x_1$ is compatible.}\label{fig:issues}
\end{figure}
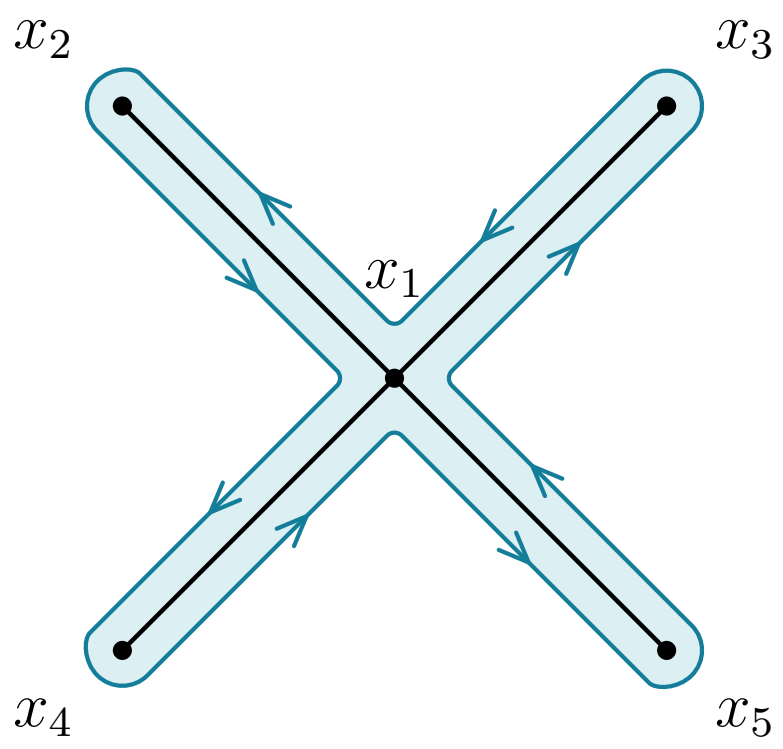

\begin{proof}
Since the balls $B_{h/8}(x_i)$ cover $\mathcal M$, choose
a spanning tree $T$ as in Lemma~\ref{lemma:embedded manifold}
and root it at $x_1$. Its edges have length at most $h/4$.

We first consider the case $d\geq 3$. 
Parametrize the doubled-tree traversal as a $\mathcal C^1$
map $\Gamma:\T\to\mathcal M$, using edge parametrizations
with velocity vanishing at their endpoints.

Assign both traversals of each edge to its child, except
that the initial half of the first outgoing traversal is
assigned to the root. This partitions $[0,1]$, up to shared
endpoints, into nonempty sets $I_1,\ldots,I_N$, each a finite
union of closed intervals, such that
\[
  \Gamma(I_i)\subseteq B_{h/4}(x_i),
  \qquad
  a_i:=\int_{I_i}\|\dot\Gamma(s)\|\,ds\leq\frac h2.
\]

Fix $0<\varepsilon<h/8$. By
Lemma~\ref{lemma:hirsch}, choose a $\mathcal C^1$ embedding
$g:\T\to\mathcal M$ satisfying
\[
  \sup_{s\in\T} d_{\mathcal M}\big(g(s),\Gamma(s)\big)<\varepsilon,
  \qquad
  \sup_{s\in\T}
  \left|\|\dot g(s)\|-\|\dot\Gamma(s)\|\right|
  <\varepsilon.
\]
Since $\varepsilon<\inj(\M)$ and $\Gamma$ is contractible, $g$ is also contractible \cite[6-16(b)]{Lee}.	

Apply Lemma~\ref{lemma:interpolation} to $g$, retaining
all endpoints of the intervals forming the $I_i$.
Choose the interpolation parameters
$0=s_0<\cdots<s_n=1$ sufficiently finely that
\[
  \ell (g|_{[s_j,s_{j+1}]})<\varepsilon,\qquad j=0,\ldots,n-1.
\]
Let $\gamma$ be the resulting simple geodesic interpolation. Then $\gamma$ is a simple geodesic cycle that is contractible because $\sup_{s\in[0,1]}d_{\M}(\gamma(s),g(s))<\inj(\M)$ and $g$ is contractible \cite[6-16(b)]{Lee}.  
Assign each interpolating segment the label of its parameter
interval, and denote the corresponding segment families
by $\gamma_i$.

For $s\in[s_j,s_{j+1}]\subseteq I_i$, the minimizing
property of the interpolating segment gives
\[
  d_{\mathcal M}(\gamma(s),x_i)
  \leq d_{\mathcal M}(\gamma(s),g(s_j))
      +d_{\mathcal M}(g(s_j),\Gamma(s_j))
      +d_{\mathcal M}(\Gamma(s_j),x_i)
  <2\varepsilon+\frac h4<\frac h2.
\]
Moreover, minimizing interpolation does not increase
the assigned lengths, so
\[
  \ell(\gamma_i)
  \leq\int_{I_i}\|\dot g(s)\|\,ds
  \leq a_i+\varepsilon|I_i|
  \leq\frac h2+\varepsilon
  \leq h.
\]
Finally, reparametrize $\gamma$ by arclength, preserving
the assigned portions, their localization, and their lengths.

\medskip
\noindent\emph{Case $d=2$.} We must find a replacement for Lemma \ref{lemma:hirsch} that provided a simple curve from the embedded tree $T$. 
The idea is to take a sufficiently thin regular
neighborhood $D_\delta$ of the tree \cite{Akitaya,LeGall}. 
This is a disc
\cite[Corollary~3.27]{RourkeSanderson}, and 
its boundary is therefore a simple contractible closed
curve, see Figure \ref{fig:tubular neighborhood}. Implementing this idea requires some technical care, as illustrated in Figure~\ref{fig:issues}. For instance, the tree traversal must follow the cyclic order induced by the boundary, and 
the approximation must preserve the assigned local length bounds. The details are given in Appendix~\ref{app:boundary}. 
\end{proof}

\begin{figure}
\begin{tikzpicture}[scale=.9,
  line cap=round,
  line join=round,
  >={Straight Barb[length=1.7mm,width=1.4mm]},
  tree edge/.style={black,line width=0.85pt}
]
  \coordinate (x1) at ( 0.0, 3.4);
  \coordinate (x2) at (-2.0, 0.4);
  \coordinate (x3) at ( 2.5, 0.4);
  \coordinate (x4) at (-3.5,-2.5);
  \coordinate (x5) at (-0.5,-2.5);

  \path[
    draw=tubeBlue,
    fill=tubeFill,
    line width=0.8pt,
    postaction={decorate},
    decoration={markings,
      mark=at position 0.0672 with
        {\arrow[draw=tubeBlue]{>}},
      mark=at position 0.2132 with
        {\arrow[draw=tubeBlue]{>}},
      mark=at position 0.3325 with
        {\arrow[draw=tubeBlue]{>}},
      mark=at position 0.4645 with
        {\arrow[draw=tubeBlue]{>}},
      mark=at position 0.5596 with
        {\arrow[draw=tubeBlue]{>}},
      mark=at position 0.6889 with
        {\arrow[draw=tubeBlue]{>}},
      mark=at position 0.7959 with
        {\arrow[draw=tubeBlue]{>}},
      mark=at position 0.9106 with
        {\arrow[draw=tubeBlue]{>}}
    }
  ]
    (2.6997,0.5664)
    -- (0.1997,3.5664)
    arc[start angle=39.8056,end angle=146.3099,radius=0.26]
    -- (-2.2163,0.5442)
    arc[start angle=146.3099,end angle=152.6501,radius=0.26]
    -- (-3.7309,-2.3806)
    arc[start angle=152.6501,end angle=332.6501,radius=0.26]
    [rounded corners=2pt]
    -- (-2.0000,-0.1659)
    -- (-0.7309,-2.6194)
    [sharp corners]
    arc[start angle=-152.6501,end angle=27.3499,radius=0.26]
    [rounded corners=2pt]
    -- (-1.6986,0.3833)
    -- (0.0232,2.9661)
    -- (2.3003,0.2336)
    [sharp corners]
    arc[start angle=-140.1944,end angle=39.8056,radius=0.26]
    -- cycle;

  \draw[tree edge]
    (x1) -- (x2)
    (x1) -- (x3)
    (x2) -- (x4)
    (x2) -- (x5);

  \node at (0.7,1.65) {$\Gamma$};

  \foreach \i in {1,...,5} {
    \filldraw[fill=black,draw=white,line width=0.45pt]
      (x\i) circle[radius=2.05pt];
  }

  \node[left=10pt]  at (x1) {$x_1$};
  \node[left=10pt]  at (x2) {$x_2$};
  \node[right=10pt] at (x3) {$x_3$};
  \node[right=10pt] at (x4) {$x_4$};
  \node[left=10pt]  at (x5) {$x_5$};

  \node[text=tubeBlue] (labelD)
    at (-3.10,1.90) {$D_\delta$};

  \draw[tubeBlue,line width=0.45pt,shorten <=3pt]
    (labelD.east) -- (-1.12,1.90);

  \node[text=tubeBlue] (labelBoundary)
    at (2.60,2.20) {$\partial D_\delta$};

  \draw[tubeBlue,line width=0.45pt,shorten <=3pt]
    (labelBoundary.west) -- (1.5997,1.8864);
\end{tikzpicture}
\caption{The tree induces a thin regular neighborhood $D_\delta$, whose boundary is a simple closed curve. 
}\label{fig:tubular neighborhood}
\end{figure}
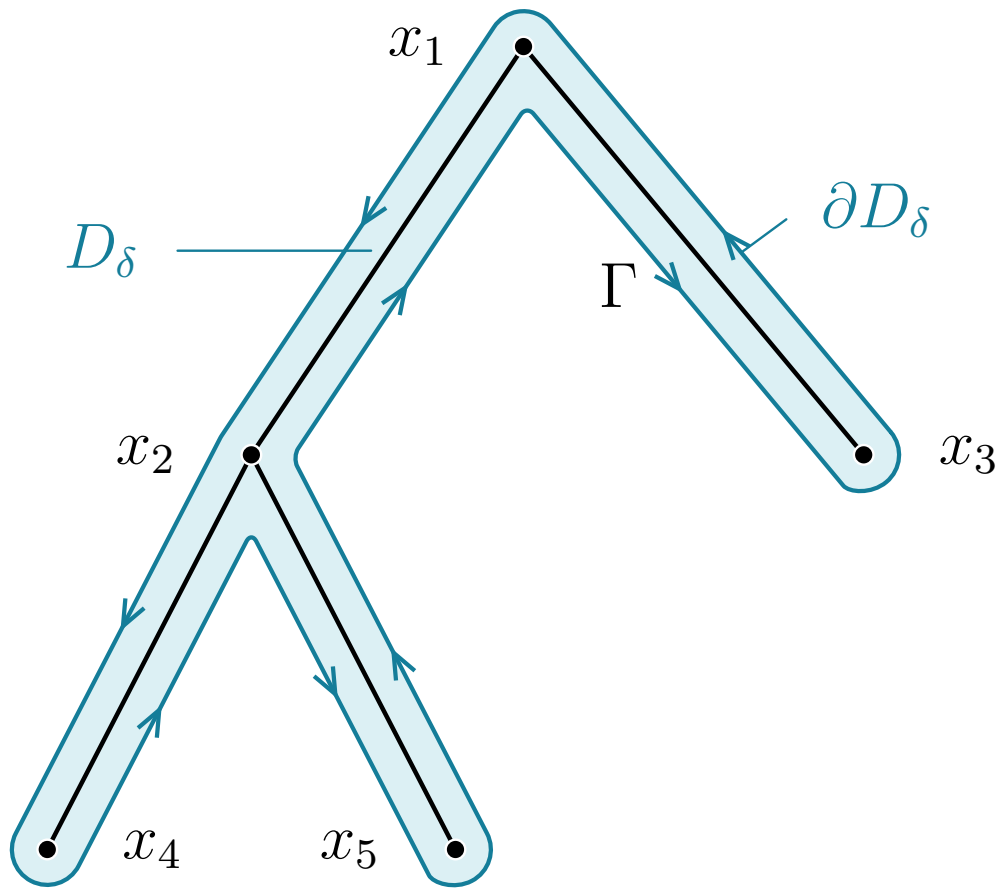

\subsection{Simple curves with adjustable local length}
Proposition \ref{prop:simple initial curve} provides a simple, contractible geodesic cycle that we will modify to meet the local length requirements. We first develop a template for local detours.

\begin{figure}
\begin{tikzpicture}[scale=1,
  line cap=round, line join=round,
    >={Straight Barb[length=2mm, width=1.7mm]}
]
  \def\H{3}       
  \def\currentH{1.5} 

  \coordinate (P)  at (0,0);
    \coordinate (T1) at (3,\H);
  \coordinate (B1) at (6,0);
    \coordinate (T2) at (7.5,\H);
  \coordinate (B2) at (9,0);
    \coordinate (T3) at (9.75,\currentH);
        \coordinate (B3) at (10.5,0);
  \coordinate (Q)  at (12,0);

  
   \node at (3,\H+0.3) {$\Phi(m_1,A)$};
  \node at (7.5,\H+0.3) {$\Phi(m_2,A)$};
    \node at (9.75,\currentH+0.3) {$\Phi(m_3,a_3)$};
  \draw[sphBlue, line width=1.1pt]
     (P) -- (T1) -- (B1)
        -- (T2) -- (B2) -- (T3) -- (B3) -- (Q) ;

  \foreach \P in {T1,T2,B1,B2}
    \fill[sphBlue] (\P) circle[radius=1.3pt];

  \foreach \P/\lab in {P/0,Q/r} {
    \fill[black!90] (\P) circle[radius=1.65pt];
    \node[below=5pt] at (\P) {$\lab$};
  }   
      \draw[<->, line width=0.65pt]
    (9.75,0.1) -- (9.75,{\currentH-0.2})
    ;

  \node at (9,-0.42) {$s_2$};
  \node at (6,-0.42) {$s_1$};
  \node at (10.5,-0.42) {$s_3$};

  \draw[dashed]
    (0,0.00)
                - - (12,0.00) 
                ;
\end{tikzpicture}
\caption{For $\nu(s)\in T_{\sigma(s)}\M$ normal to $\dot{\sigma}(s)$, set $\Phi(s,a):=\exp_{\sigma(s)}(a\nu(s))$. 
For $d_{\M}(p,q)=r$ and $\vartheta=1/2$, we illustrate two teeth of full height and one tooth grown to half height.}\label{fig:m}
\end{figure}
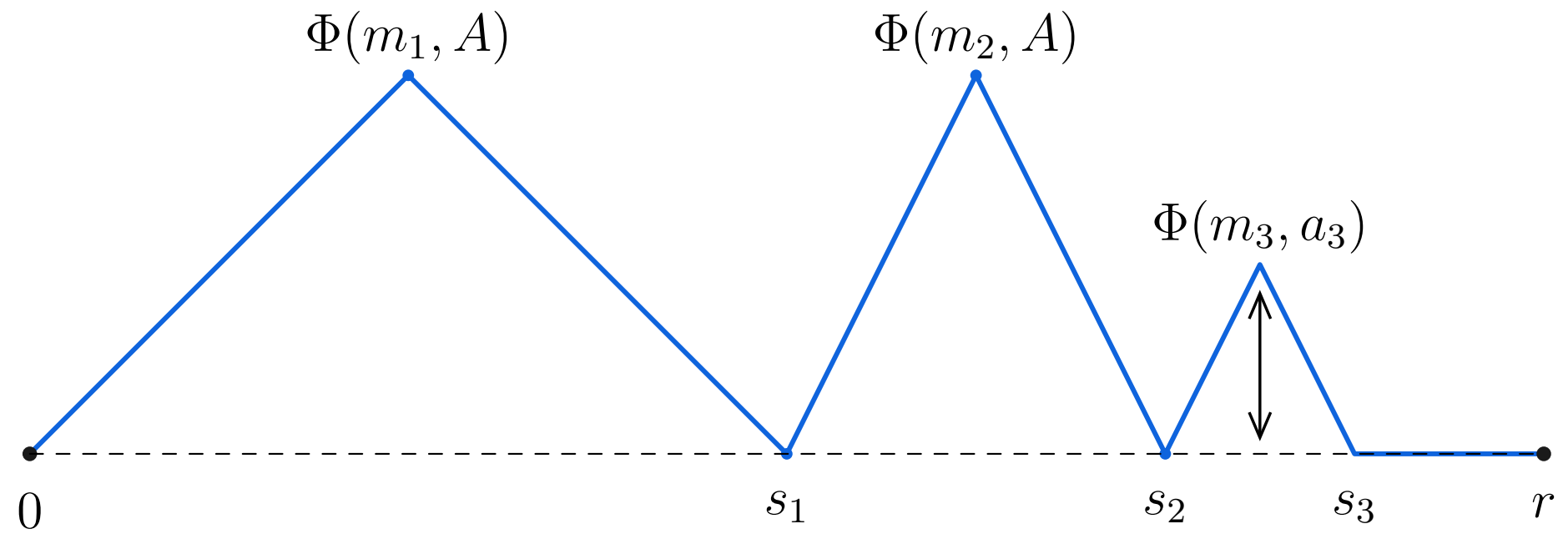

Related local perturbations of embedded spherical arcs with
continuously prescribed lengths appear in
\cite[Lemma~3.1]{GhomiRaffaelli}. For background on homotopies of paths with fixed endpoints,
see \cite[Chapter~7]{Lee2}.
\begin{proposition}[Local detours of prescribed length]
\label{prop:local-detours}
Let $\mathcal M$ be a closed, connected Riemannian
manifold of dimension $d\geq2$. Let $p,q\in\mathcal M$ satisfy
\[
  0<d_{\mathcal M}(p,q)=:r<\operatorname{inj}(\mathcal M),
\]
and let $\sigma:[0,r]\to\mathcal M$ be the unit-speed
minimizing geodesic from $p$ to $q$.

For every open neighborhood $U$ of $\sigma([0,r])$, there
exists a continuous map
\[
  \alpha:[r,\infty)\times[0,1]\longrightarrow U
\]
such that, for every $L\in[r,\infty)$, the curve
$\alpha_L:=\alpha(L,\cdot)$ is simple, piecewise geodesic, homotopic
to $\sigma$, 
and parametrized at constant speed $L$, with
\[
  \alpha_L(0)=p,\qquad
  \alpha_L(1)=q,\qquad
  \ell(\alpha_L)=L.
\]
Moreover, the map $[r,\infty)\ni L
  \longmapsto
  \int_{\alpha_L}f$ 
is continuous, for all $f\in\mathcal C(\M)$.
\end{proposition}
\begin{proof}
We construct the curves by successively growing geodesic teeth
along $\sigma$ until the prescribed length is attained, cf.~Figure \ref{fig:m}.

Along $\sigma$, choose a smooth parallel unit normal field
$\nu(s)\in T_{\sigma(s)}\M$, and set
\[
  \Phi(s,a):=\exp_{\sigma(s)}(a\nu(s)).
\]
After slightly extending $\sigma$ past its endpoints,
the normal tubular-neighborhood theorem and the normal
distance formula
\cite[Theorem~5.25 and Corollary~6.39(d)]{Lee}
provide a normal tube $V\subseteq U$, a smooth longitudinal
coordinate $S$ on $V$, and a sufficiently small $A>0$ such that
\begin{equation}\label{eq:dM etc}
  S(\Phi(s,a))=s,\qquad
  d_\M(\Phi(s,a),\sigma([0,r]))=a,
  \qquad s\in[0,r],\quad a\in[0,A].
\end{equation}

Choose $A>0$ small enough that $r+A<\inj(\M)$, and denote
the constant-speed minimizing geodesic from $\sigma(s)$
to $\Phi(m,a)$ by $c_{s,m,a}:[0,1]\to\M$. To verify that $c_{s,m,a}(u)\in V$, we need some preparation. 
The segments $c_{s,m,a}$ depend smoothly on $(s,m,a,u)$
\cite[Problem~6-16(a)]{Lee}. At $a=0$, we have
\[
  c_{s,m,0}(u)=\sigma((1-u)s+u m).
\]
Since $\sigma([0,r])$ is a compact subset of the open
tube $V$, uniform continuity ensures that, after decreasing
$A$, all segments $c_{s,m,a}$ lie in $V$, uniformly for $s,m\in[0,r]$, $a\in[0,A]$, and $u\in[0,1]$.

Since $c_{m,m,a}(u)=\Phi(m,ua)$, we have
$S(c_{m,m,a}(u))=m$. Thus the smooth function
\[
F(s,m,a,u):=\partial_u S(c_{s,m,a}(u))
\]
vanishes when $s=m$. The identity
\[
\frac{F(s,m,a,u)}{m-s}
=
-\int_0^1
\partial_1F\bigl(m+t(s-m),m,a,u\bigr)\,dt,
\qquad s\ne m,
\]
therefore gives a smooth extension of the quotient across
$s=m$. At $a=0$ it equals $1$.
Decreasing $A$ again makes it positive uniformly.
Thus, the map $u\mapsto S(c_{s,m,a}(u))$ is strictly
increasing when $s<m$ and strictly decreasing when $s>m$.

For $0=s_0<s_1<s_2<\cdots <r$ and $m_j:=\frac{s_{j-1}+s_j}{2}$, replace $\sigma|_{[s_{j-1},s_j]}$ by joining consecutive vertices 
\[
  \sigma(s_{j-1})\rightarrow\Phi(m_j,a)
  \rightarrow\sigma(s_j),
\]
by minimizing geodesics. We refer to this geodesic chain as the $j$-th tooth. Its length
\[
  T_j(a):=
  d_\M(\sigma(s_{j-1}),\Phi(m_j,a))
  +d_\M(\Phi(m_j,a),\sigma(s_j))
\]
is continuous with $T_j(0)=s_j-s_{j-1}$, and \eqref{eq:dM etc} implies $T_j(a)\geq2a$. 

Consider the smooth squared-distance function
\[
  D(s,m,a):=d_\M(\sigma(s),\Phi(m,a))^2.
\]
Since $D(m,m,a)=a^2$, smoothness and compactness give
$\delta>0$ such that
\[
  \partial_a^2D(s,m,a)\geq1,
  \qquad
  |s-m|<\delta,\quad 0\leq a\leq A.
\]
Choose $\vartheta\in(0,1)$ sufficiently close to $1$ that
$r(1-\vartheta)<2\delta$, and specify $s_j$ according to 
\[
  s_j:=r(1-\vartheta^j),
\]
see Figure \ref{fig:m} for $\vartheta=1/2$. 

Then $|s_{j-1}-m_j|=|s_j-m_j|<\delta$ for every $j$.
Orthogonality and the first variation formula
\cite[Theorem~6.3]{Lee} give $\partial_aD(s,m,0)=2(m-s)\langle \dot{\sigma}(m),\nu(m)\rangle = 0$.
Hence, 
\[
  \partial_aD(s,m,a)
  =\int_0^a\partial_a^2D(s,m,b)\,db
  \geq a>0,\qquad 0<a\leq A,
\]
for either endpoint $s$ of each tooth.
Consequently, both sidelengths, and therefore $T_j(a)$,
are strictly increasing on $[0,A]$ for every $j$. 

The longitudinal coordinate increases strictly along the
oriented tooth, from $s_{j-1}$ through $m_j$ to $s_j$.
Thus each tooth is simple, distinct teeth meet only at
shared base endpoints, and there are no intersections
with the retained portions of $\sigma$.

During $\tau\in[n,n+1]$, grow the $(n+1)$-th tooth from height
zero to $A$, keeping earlier teeth at fixed $A$. For $\tau\in[0,\infty)$, define the height of the $j$-th tooth by 
\[
  a_j(\tau)=A\min\big(1,\max(0,\tau-j+1)\big).
\]
For fixed $\tau$, only finitely many teeth have nonzero height and $a_j(\tau)=0$ for all $j\geq \tau+1$.
Since $s_j-s_{j-1}=r\vartheta^{j-1}(1-\vartheta)$, the total length
\[
  \lambda(\tau)
  =\sum_{j=1}^n T_j(A) \, + T_{n+1}(a_{n+1}(\tau)) +\vartheta^{n+1}r,\qquad \tau\in[n,n+1].
\]
is continuous and strictly increasing with 
$\lambda(0)=r$ and $\lambda(n)\geq 2nA\to\infty$ for $n\to\infty$. Thus $\lambda:[0,\infty)\to[r,\infty)$ is a homeomorphism.

 Let $\beta_\tau:[0,1]\to U$ be the resulting chain
parametrized at constant speed $\lambda(\tau)$. On every bounded $\tau$-interval, only finitely many teeth
are involved. By including the corresponding base endpoints
and midpoint vertices also before a tooth starts growing,
the chains have a fixed finite number of geodesic pieces
with continuously varying endpoints and positive lengths.
Their constant-speed parametrizations $\beta_\tau$ therefore
depend jointly continuously on $(\tau,s)$.

Define
\[
  \alpha_L(s):=\beta_{\lambda^{-1}(L)}(s),
  \qquad L\in[r,\infty),\quad s\in [0,1].
\]
Then $\alpha_L$ is simple, piecewise
geodesic, has endpoints $p,q$, and has constant speed
and length $L$. Moreover, $\alpha_r(s)=\sigma(rs)$ and,
for every $f\in\mathcal C(\M)$,
\[
  \int_{\alpha_L}f
  =L\int_0^1f(\alpha_L(s))\,ds
\]
depends continuously on $L$.
Finally,
\[
  H_L(\theta,s):=\alpha_{r+\theta(L-r)}(s),\qquad \theta,s\in[0,1],
\]
is a homotopy in $U$ from $s\mapsto \sigma(rs)$
to $\alpha_L$.
\end{proof}

We obtain the desired curves with locally adjustable lengths by
starting from the initial curve in Proposition~\ref{prop:simple initial curve}
and adding the local detours provided by Proposition~\ref{prop:local-detours}.

\begin{theorem}\label{thm:curves are done}
Let $R_1,\dots,R_N$ be a partition with points $x_1,\ldots,x_N$ as in Theorem \ref{thm:partition} and suppose that $h=8bN^{-1/d}<\frac{2}{3}\inj(\M)$. Then there is a family of simple, 
contractible, geodesic cycles 
$\{\gamma_w\}_{w\in[1,3]^N}$, whose segment traversals can be
partitioned into families $\gamma_{i,w}$ satisfying
\[
  \gamma_{i,w}\subseteq B_{h}(x_i),
  \qquad \ell(\gamma_{i,w})=hw_i.
\]
Moreover, the map $[1,3]^N\ni 
  w\longmapsto
  \int_{\gamma_w}f$ 
is continuous, for all $f\in\mathcal C(\M)$.
\end{theorem}
\begin{proof}
Let $\gamma$ be as in Proposition \ref{prop:simple initial curve},
with assigned segment families $\gamma_i$ satisfying
\[
\gamma_i\subseteq B_{h/2}(x_i),
\qquad
\ell(\gamma_i)\leq h.
\]
The construction ensures that each $\gamma_i$ contains a
nondegenerate geodesic segment. For each $i$, choose a closed
subsegment $\sigma_i$ strictly inside such a segment, with the
orientation inherited from $\gamma$, and put
\[
r_i:=\ell(\sigma_i)<\operatorname{inj}(M).
\]
Choose pairwise disjoint open tubular neighborhoods
$U_i\subseteq B_h(x_i)$ of $\sigma_i$, sufficiently small that
each $U_i$ meets $\gamma$ only along the original geodesic
segment containing $\sigma_i$.

Apply the construction in the proof of
Proposition \ref{prop:local-detours} to $\sigma_i$ inside $U_i$,
and denote the resulting detours by $\alpha_{i,L}$, $L\geq r_i$.
For $w\in[1,3]^N$, define
\[
L_i(w_i):=r_i+hw_i-\ell(\gamma_i).
\]
Since $w_i\geq1$ and $\ell(\gamma_i)\leq h$, we have
$L_i(w_i)\geq r_i$. Replace each $\sigma_i$ by
$\alpha_{i,L_i(w_i)}$, and denote the resulting cycle by
$\gamma_w$ and its assigned segment families by $\gamma_{i,w}$.
Then
\[
\ell(\gamma_{i,w})
=\ell(\gamma_i)-r_i+L_i(w_i)
=hw_i.
\]
Moreover, $\gamma_{i,w}\subseteq B_h(x_i)$, since the unchanged
segments lie in $B_{h/2}(x_i)$ and the detour lies in $U_i$.

In the construction of Proposition \ref{prop:local-detours},
the longitudinal coordinate increases strictly along each
detour between its attachment points. The retained portions
of the original geodesic segment lie beyond these endpoint
levels. Thus each detour meets the retained curve only at its
endpoints. Together with the disjointness of the neighborhoods
$U_i$, this shows that $\gamma_w$ is simple. Each replacement
is homotopic to $\sigma_i$ with endpoints fixed, so $\gamma_w$
is homotopic to $\gamma$ and is therefore contractible.

For every $f\in \mathcal C(\M)$, we have
\[
\int_{\gamma_w}f
=
\int_\gamma f
-\sum_{i=1}^N\int_{\sigma_i}f
+\sum_{i=1}^N\int_{\alpha_{i,L_i(w_i)}}f.
\]
Each term in the last sum depends continuously on $w_i$ by
Proposition \ref{prop:local-detours}. Hence
$w\mapsto\int_{\gamma_w}f$ is continuous.
Finally, reparametrize each $\gamma_w$ by arclength.
\end{proof}

This theorem accomplishes Step (i) in the structural Theorem \ref{thm:foundation}. Step (ii) asks for local oscillation estimates.

\section{Local oscillation estimates and proofs of the main theorems}\label{sec:osc manifold}
We verify condition (ii) of Theorem \ref{thm:foundation} in order to prove Theorem \ref{thm:manifold-design-curves} in the introduction. 
\subsection{Average local oscillation estimates}\label{sec:osc}
We use the $L^1$-Marcinkiewicz-Zygmund estimates from \cite[Theorem 5.1]{FilbirMhaskar} as used in \cite{Gariboldi1}.
\begin{lemma}[\cite{FilbirMhaskar}]
There is a constant $A_b\geq 1$ such that, for every partition $R_1,\dots,R_N$ with points $x_1,\ldots,x_N$ as in Theorem \ref{thm:partition} and for all positive real $t\leq N^{1/d}$,
\begin{equation}\label{eq:manifold-mz}
  \left|\frac1N\sum_{i=1}^N |P(y_i)|-\|P\|_{L^1(\mathcal M)}\right|
  \leq A_b tN^{-1/d}\|P\|_{L^1(\mathcal M)},
  \quad  \forall \, y_i\in R_i,\quad P\in\Pi_t(\mathcal M).
\end{equation}
\end{lemma}
Although these estimates were also used in \cite{Gariboldi1}, we take a different route, obtaining bounds on the local oscillations directly.
\begin{proposition}\label{lem:manifold-oscillation}
For a partition $R_1,\dots,R_N$ with points $x_1,\ldots,x_N$ as in Theorem \ref{thm:partition}, put $V_i:=B_{8bN^{-1/d}}(x_i)$, so that $R_i\subseteq V_i$. If $t\leq N^{1/d}$, then
\begin{equation}\label{eq:manifold-oscillation}
  \frac1N\sum_{i=1}^N\osc_{V_i}(P)
  \leq 4A_{8b}\,t N^{-1/d} \|P\|_{L^1(\mathcal M)},
  \qquad P\in\Pi_t(\M).
\end{equation}
\end{proposition}
\begin{proof}
Choose $u_i,v_i\in V_i$ where $|P|$ attains its maximum and
minimum on $V_i$, respectively. The sets
$\tilde{R}_i=R_i\cup\{u_i,v_i\}\subseteq V_i$ still form an equal-measure partition with $\tilde{R}_i\subseteq B_{8bN^{-1/d}}(x_i)$.  
Apply \eqref{eq:manifold-mz} to the choices $u_i$ and $v_i$ and
subtract the resulting bounds:
\[
  \frac1N\sum_{i=1}^N\osc_{V_i}(|P|)
  =\frac1N\sum_{i=1}^N\bigl(|P(u_i)|-|P(v_i)|\bigr)
  \leq2A_{8b}\, tN^{-1/d}\|P\|_{L^1(\mathcal M)}.
\]
If $P$ has a fixed sign
on $V_i$, its oscillation equals that of $|P|$. Otherwise,
$P$ vanishes somewhere in $V_i$, and
\begin{equation*}
\osc_{V_i}(P)\leq2\max_{V_i}(|P|)
=2\osc_{V_i}(|P|).
\end{equation*}
Thus $\osc_{V_i}(P)\leq
2\osc_{V_i}(|P|)$ in both cases, proving the claim.
\end{proof}

The oscillation estimates are used in combination with suitable parameter choices to verify the requirement $\frac{12}{N}\sum_{i=1}^N \osc_{V_i}(P) < 1$ in \eqref{eq:osc finale}.

\subsection{Parameter choice and proof of Theorem \ref{thm:manifold-design-curves}}
We verify Theorem \ref{thm:manifold-design-curves} by applying the structural Theorem \ref{thm:foundation}. This still requires a suitable choice of the parameters $N$ and $h$. 

\begin{lemma}[Choice of parameters]\label{lemma:parameters manifold}
Suppose $d\geq 2$. Fix the constants $b,A_{8b}$
above and choose $C_{\mathcal M}
    :=16b\max\big(
        64A_{8b},\frac{16b}{\inj(\M)}
    \big)^{d-1}$. For $t\geq1$ and $L\geq C_{\mathcal M}\,t^{d-1}$, define
\[
  \qquad N:=\Big\lceil \Big(\frac{L}{16b}\Big)^{d/(d-1)}\,\Big\rceil,
  \qquad h:=8bN^{-1/d},
  \]
Then we have $h<\frac{2}{3}\inj(\M)$ and 
\begin{equation}\label{eq:manifold-parameters}
t N^{-1/d} \leq\frac 1{64A_{8b}},\qquad \frac{L}{hN} \in \Big[\frac{3}{2},\frac{5}{2}\Big].
\end{equation}
\end{lemma}
\begin{proof}
The lower bound on $L$ gives
\[
    h=8bN^{-1/d}
    \le 8b\Big(\frac{16b}{L}\Big)^{1/(d-1)}
    \le \frac{\inj(\M)}{2t}
    \le \frac12\inj(\M).
\]

Set $\tilde{N}:=\left(\frac{L}{16b}\right)^{d/(d-1)}\ge(64A_{8b}\,t)^d$.
Since $\tilde{N}\leq N$, we derive
\[
tN^{-1/d}\le t\tilde{N}^{-1/d}\le\frac1{64A_{8b}}.
\]
Using \(hN=8bN^{(d-1)/d}\), we obtain
\[
\frac{L}{hN}
=2\Big(\frac{\tilde{N}}{N}\Big)^{(d-1)/d}
\in
\Big[2\left(\frac34\right)^{(d-1)/d},\,2\Big]
\subseteq
\left[\frac32,\frac52\right].
\]
\end{proof}

We have all ingredients to apply Theorem \ref{thm:foundation}.
\begin{proof}[Proof of Theorem \ref{thm:manifold-design-curves}]
Choose $C_\M$ as in Lemma \ref{lemma:parameters manifold} and for $t\geq1$ and $L\geq C_{\mathcal M}\,t^{d-1}$ choose also $N$ and $h$ as in Lemma \ref{lemma:parameters manifold}, so that $\frac{L}{hN}\in [\frac{3}{2},\frac{5}{2}]$ and $t N^{-1/d} \leq\frac 1{64A_{8b}}$ hold. 

Theorem \ref{thm:curves are done} provides a family of simple, contractible geodesic cycles $\{\gamma_w\}_{w\in[1,3]^N}$, each $\gamma_w$ partitioned into $\gamma_{i,w}$  satisfying
\[
  \gamma_{i,w}\subseteq B_{h}(x_i),
  \qquad \ell(\gamma_{i,w})=hw_i,
\]
as well as the continuity conditions. 
For $V_i=B_{8bN^{-1/d}}(x_i)$, we have $ \gamma_{i,w}\subseteq V_i$ and $R_i\subseteq V_i$. Proposition \ref{lem:manifold-oscillation}, for $P\in\Pi_t(\M)$ with $\|P\|_{L^1(\M)}=1$, yields
\begin{equation*}
\frac{12}{N}\sum_{i=1}^N\osc_{V_i}(P) \leq 48\,A_{8b}\, t N^{-1/d} \leq 3/4<1.
\end{equation*}
Thus, all assumptions of the structural Theorem \ref{thm:foundation} are satisfied, and there is $\hat{w}\in[1,3]^N$ such that 
$\gamma_{\hat{w}}$ is a $t$-design curve of length $L$. This proves Theorem \ref{thm:manifold-design-curves} in the introduction.
\end{proof}

\subsection{Smooth $t$-design curves}
We prove Theorem~\ref{thm:manifold-design-curves2}, which establishes the existence of simple, smooth, contractible $t$-design curves.

The idea is to use the family of simple, contractible geodesic cycles $(\gamma_w)_{w\in[1,3]^N}$ from Theorem \ref{thm:curves are done} that satisfy the positivity condition for Brouwer's degree theorem. By corner smoothing introducing only sufficiently small perturbations $\varepsilon>0$, we derive a family of simple, smooth, contractible curves $(\gamma_{w,\varepsilon})_{w\in[1,3]^N}$. Sufficient control over the smoothing process ensures that this smooth family still satisfies the positivity condition and each curve is still simple and contractible.

The following lemma is used for $W=[1,3]^N$.
\begin{lemma}[Parametric smoothing of geodesic cycles]
\label{lem:parametric-smoothing}
Let $W$ be a compact metric space, and let
\[
    p_0,\ldots,p_{m-1}:W\longrightarrow\mathcal M,
\]
be continuous maps, with $p_m=p_0$, such that $0<
    d_{\mathcal M}\bigl(p_j(w),p_{j+1}(w)\bigr)
    <\operatorname{inj}(\mathcal M)$, for $w\in W$ and $0\leq j<m$.
For $w\in W$, let
$\gamma_w:\mathbb R/m\mathbb Z\to\mathcal M$ be the geodesic
cycle obtained by joining $p_j(w)$ to $p_{j+1}(w)$ by the
minimizing geodesic, parametrized with constant speed
on $[j,j+1]$, so that $\gamma_w(j)=p_j(w)$.
Suppose that every $\gamma_w$ is simple.

Then, for every $\varepsilon>0$, there exists a continuous map
\begin{equation}\label{eq:WCk}
    W\longrightarrow
    \mathcal C^1(\mathbb R/m\mathbb Z,\mathcal M),
    \qquad
    w\longmapsto\gamma_{w,\varepsilon},
\end{equation}
with the following properties:
\begin{enumerate}
    \item
    Each $\gamma_{w,\varepsilon}$ is a smooth embedding, homotopic to $\gamma_w$, and 
    \[
        d_{\mathcal M}\bigl(
            \gamma_{w,\varepsilon}(s),\gamma_w(s)
        \bigr)
        \leq \varepsilon,\qquad \forall\, s\in\mathbb R/m\mathbb Z,\quad w\in W.
    \]
    \item
    The path integrals\,\footnote{For $\gamma_w:[0,m]\to\M$, the path integral is $\int_{\gamma_w} P = \int_0^m P(\gamma_w(s)) |\dot{\gamma}_w(s)|\dd s$.} 
 satisfy
    \[
        \Big|
            \int_{\gamma_{w,\varepsilon}} P
            -
            \int_{\gamma_w} P
        \Big|
        \leq \varepsilon\|P\|_\infty,
        \qquad
        \forall \,P\in \mathcal C(\mathcal M),\quad w\in W.
    \]
\end{enumerate}
\end{lemma}

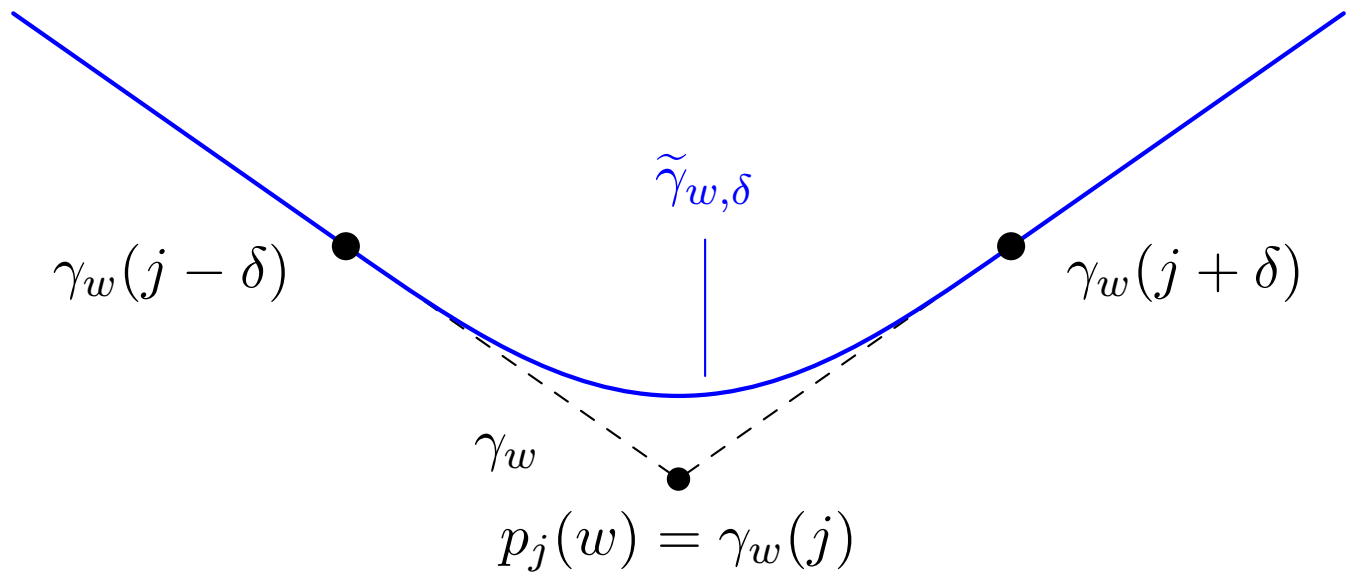
\begin{figure}[htbp]
\centering
\begin{tikzpicture}[scale=2.2,
    line cap=round,
    line join=round,
    font=\small
]
    \draw[blue, thick]
        (-2,1.4) -- (-1,.7);
    \draw[blue, thick]
        (1,.7) -- (2,1.4);

    \draw[black, dashed]
        (-1,.7) -- (0,0) -- (1,.7);

    \draw[blue, thick]
        (-1,.7)
        .. controls (-.65,.455) and (-.35,.25) .. (0,.25)
        .. controls (.35,.25) and (.65,.455) .. (1,.7);

    \fill[black] (-1,.7) circle (1.2pt);
    \fill[black] (1,.7) circle (1.2pt);

    \node[anchor=east] at (-1.1,.62)
        {$\gamma_w(j-\delta)$};
    \node[anchor=west] at (1.1,.62)
        {$\gamma_w(j+\delta)$};

    \fill[black] (0,0) circle (1pt);
    \node[below] at (0,-.03) {$p_j(w)=\gamma_w(j)$};
    \node[black, anchor=east] at (-.35,.08) {$\gamma_w$};

    \node[blue] at (.08,.89)
        {$\widetilde{\gamma}_{w,\delta}$};
    \draw[blue, thin] (.08,.72) -- (.08,.31);
\end{tikzpicture}
\caption{The corner pieces from $\gamma_w(j-\delta)$ to $p_j(w)$ and from $p_j(w)$ to $\gamma_w(j+\delta)$ are replaced by the smooth blue arc. 
The curves $\widetilde{\gamma}_{w,\delta}$ and $\gamma_w$ agree outside the parameter intervals
$[j-\delta,j+\delta]$.}
\label{fig:corner-smoothing}
\end{figure}

The proof of Lemma \ref{lem:parametric-smoothing} is based on smoothing of corners by convolution with a mollifier in normal coordinates as in \cite[Section~3]{EpsteinONeil2016}, see also \cite[Chapter~2, Theorems~2.3 and 2.5]{Hirsch}. We provide a sketch of the proof and refer to Appendix \ref{app:corner} for the complete arguments, in particular, for the choice of $\varepsilon$ independent of $w$.
\subsubsection*{Sketch of the proof}
Fix a nonnegative, even function $\eta\in \mathcal C_c^\infty((-1,1))$ with $\int_{-1}^1\eta(s)\dd s=1$ and put $\eta_\delta(s)=\frac{1}{\delta}\eta(\frac{s}{\delta})$ for $0<\delta<1/2$. Let $v_{j,+}(w)$ and $v_{j,-}(w)$ denote
the outgoing and incoming velocities of $\gamma_w$ at $p_j(w)$ in normal coordinates. Write the $j$-th corner piece in these normal coordinates centered
at $p_j(w)$ as $q_{j,w}(s)$, and extend it to 
\[
q_{j,w}(s)=
\begin{cases}
(s-j)v_{j,-}(w),&s\le j,\\
(s-j)v_{j,+}(w),&s\ge j.
\end{cases}
\]
The smoothed curve in these normal coordinates is defined as $q_{j,w,\delta}:=q_{j,w}*\eta_\delta$.  Its restriction to \([j-\delta,j+\delta]\) is used for the replacement of the corner piece. Since the mollifier is even with unit integral, we have 
$q_{j,w,\delta}=q_{j,w}$ for $|s-j|\geq\delta$. 

Applying this construction at every corner and leaving the
remaining geodesic pieces unchanged defines
$\widetilde\gamma_{w,\delta}$, which agrees with $\gamma_w$
outside the union of the intervals
$[j-\delta,j+\delta]$, $j=0,\ldots,m-1$, see Figure \ref{fig:corner-smoothing}. The choice of $\varepsilon$ in relation to $\delta$ and its choice independent of $w$ are provided in the detailed proof in Appendix \ref{app:corner}.

\begin{proof}[Proof of Theorem \ref{thm:manifold-design-curves2}]
Choose $C_{\mathcal M}$ as in Lemma \ref{lemma:parameters manifold}, fix $t\geq1$ and
$L\geq C_{\mathcal M}t^{d-1}$, and choose $N,h$ accordingly.
Take the family from Theorem \ref{thm:curves are done}. On the bounded parameter
ranges used here, the construction in Proposition \ref{prop:local-detours} has
a fixed finite number of continuously varying vertices.
After a fixed subdivision and reparametrization, Lemma \ref{lem:parametric-smoothing}
therefore applies with $W=[1,3]^N$.

Use $\Omega$, $w(P)$, and $F(P)$ from Section~2 with
$\Pi(\mathcal M)=\Pi_t(\mathcal M)$. 
The final estimate in the proof of Theorem \ref{thm:foundation} and the
oscillation bound in Section \ref{sec:osc} give
\[
    \langle P,F(P)\rangle_{L^2(\mathcal M)}
    \geq \frac{hN}{16},
    \qquad P\in\partial\Omega.
\]
Since $\Pi_t(\mathcal M)$ is finite-dimensional, there is
$C_t>0$ such that
$\|P\|_\infty\leq C_t\|P\|_{L^1(\mathcal M)}$
for all $P\in\Pi_t(\mathcal M)$.
Choose $\varepsilon>0$ with $\varepsilon C_t<hN/16$
and apply Lemma \ref{lem:parametric-smoothing} to obtain the smoothed family $(\gamma_{w,\varepsilon})_{w\in[1,3]^N}$, so that $\gamma_{w,\varepsilon}:[0,m]\to\M$. 

Let $K_t$ be the reproducing kernel of $\Pi_t(\M)$, and $w:\Omega\to[1,3]^N$ as in \eqref{eq:wi}. 
Define
\begin{equation*}
F_\varepsilon:\Omega\to \Pi(\M),\qquad P\mapsto \int_{0}^{m} K_t(\cdot,\gamma_{w(P),\varepsilon}(s)) |\dot{\gamma}_{w(P),\varepsilon}(s)|\dd s - L,
\end{equation*}
so that the continuity of \eqref{eq:WCk} supplied by Lemma \ref{lem:parametric-smoothing} implies that
$F_\varepsilon$ is continuous. The respective identity \eqref{eq:satisfies this here} is 
\[
    \langle Q,F_\varepsilon(P)\rangle_{L^2(\mathcal M)}
    =
    \int_{\gamma_{w(P),\varepsilon}}Q
    -
    L\int_{\mathcal M}Q\,d\mu,
    \qquad Q\in\Pi_t(\mathcal M).
\]
By part~(2) of Lemma \ref{lem:parametric-smoothing},
\[
\begin{aligned}
    \langle P,F_\varepsilon(P)\rangle_{L^2(\mathcal M)}
    &\geq
    \langle P,F(P)\rangle_{L^2(\mathcal M)}
    -\varepsilon\|P\|_\infty\\
    &\geq \frac{hN}{16}-\varepsilon C_t>0,
    \qquad P\in\partial\Omega.
\end{aligned}
\]
Theorem \ref{thm:brouwer} therefore gives $P_0\in\Omega$ with
$F_\varepsilon(P_0)=0$. Hence, for
$\gamma:=\gamma_{w(P_0),\varepsilon}$,
\[
    \int_\gamma Q
    =
    L\int_{\mathcal M}Q\,d\mu,
    \qquad Q\in\Pi_t(\mathcal M).
\]
Taking $Q\equiv1$ yields $\ell(\gamma)=L$.
By Part (1) of Lemma \ref{lem:parametric-smoothing}, $\gamma$ is a smooth embedding homotopic
to the contractible curve $\gamma_{w(P_0)}$.
Reparametrization by arclength completes the proof of Theorem \ref{thm:manifold-design-curves2}.
\end{proof}

\subsection{From spectral to spherical designs}\label{app:M as T and S}
The spherical case is a consequence of the manifold result in Theorem \ref{thm:manifold-design-curves}, whose proof we just completed. The sphere $\Sd$ is a closed, connected Riemannian manifold, but to derive Theorem \ref{thm:piecewise}, we must still relate the polynomials of total degree $\leq t$ with the spectral space $\Pi_t(\mathbb S^d)$. This is standard in the harmonic analysis of the sphere \cite{SteinWeiss}, and we recall the relevant facts. 

The space $\mathcal H_l$ of real spherical harmonics of degree $l$
is an eigenspace of the nonnegative Laplace-Beltrami operator,
\[
  \Delta_{\mathbb S^d}Y=l(l+d-1)Y,
  \qquad Y\in\mathcal H_l.
\]
Thus, for every positive integer $t$, we have 
\begin{equation*}
\Pi_t(\mathbb S^d)
  =\bigoplus_{\substack{l \geq0\\ l(l+d-1)\leq t^2}}
      \mathcal H_l
      \end{equation*}
      and may define 
\[
\mathcal P_t(\mathbb S^d)
  :=\bigoplus_{l=0}^t\mathcal H_l.
\]
The space $\mathcal P_t(\mathbb S^d)$ is known to coincide with the polynomials of total degree $\leq t$, restricted to the sphere $\Sd$. For $t\geq1$ and $d\geq 2$, the relations $t\leq\sqrt{t(t+d-1)}\leq\sqrt d\,t$ lead to
\begin{equation}\label{eq:Pi and Pi}
\Pi_t(\Sd)  \subseteq \mathcal P_{t} (\Sd)\subseteq \Pi_{\sqrt{d}\, t} (\Sd).
\end{equation}
Consequently, a spherical $t$-design curve in the sense of
Definition \ref{def:designs S}, that is with respect to $\mathcal P_{t} (\Sd)$, is a spectral $t$-design curve in the sense of Definition \ref{def:manifold-design-curve}, that is with respect to $\Pi_t(\Sd)$. On the other hand, a spectral
$\sqrt{d}\, t$-design curve on $\Sd$ is a spherical $t$-design curve. Therefore, Theorem \ref{thm:piecewise} follows from Theorem \ref{thm:manifold-design-curves}.

\section{Marcinkiewicz-Zygmund inequalities and Logvinenko-Sereda sets}\label{sec:further}

We provide two observations: the curve family in Theorem \ref{thm:curves are done} satisfies uniform Marcinkiewicz-Zygmund inequalities, and Theorem \ref{thm:foundation} leads to Logvinenko-Sereda sets with additional exact integration properties. 

\subsection{Marcinkiewicz-Zygmund inequalities}

Marcinkiewicz-Zygmund inequalities \cite{FilbirHielscherJahnUllrich,FilbirMhaskar,Grochenig2020,MarzoMZ,MhaskarNarcowichWard} compare norms obtained by sampling functions with their $L^p$-norms with respect to the underlying measure, here the Riemannian measure $\mu$ on $\M$. For $P\in\mathcal C(\M)$, define
\[
    \|P\|_{L^p(\M,\gamma_w)}
    :=
    \begin{cases}
        \displaystyle
        \Big(
            \frac{1}{\ell(\gamma_w)}
            \int_{\gamma_w}|P|^p
        \Big)^{1/p},
        & 1\le p<\infty,\\[2ex]
        \displaystyle
        \max_{0\le s\le\ell(\gamma_w)}|P(\gamma_w(s))|,
        & p=\infty.
    \end{cases}
\]
We ask when this curve-based norm is equivalent to $\|P\|_{L^p(\M,\mu)}$ on $\Pi_t(\M)$, with constants independent of $t$ and $p$.

For the sphere $\Sd$, nonsimple curves of length $\mathcal O(t^{d-1})$ satisfying Marcinkiewicz-Zygmund inequalities were constructed in \cite{EhlerGroechenigKarner}. Here we extend the setting to closed, connected Riemannian manifolds and obtain families of simple, contractible geodesic cycles satisfying such inequalities.

\begin{proposition}[Marcinkiewicz-Zygmund inequalities for the curve family]
\label{prop:MZ-curves}
For all $\varepsilon\in(0,1)$, there is a constant $c_{\M,\varepsilon}>0$ such that for the family of simple, contractible geodesic cycles $\{\gamma_w\}_{w\in[1,3]^N}$ constructed in Theorem \ref{thm:curves are done}, we have, for all positive $ t \leq  c_{\M,\varepsilon}N^{1/d}$ and every
$1\le p\le\infty$ and all $ w\in[1,3]^N$,
\begin{equation}\label{eq:MZ-curves}
   \tfrac{1}{3}(1-\varepsilon)\|P\|_{L^p(\M)}
    \le
    \|P\|_{L^p(\M,\gamma_w)}
    \le
    3 (1+\varepsilon)\|P\|_{L^p(\M)},
    \qquad \forall\, P\in\Pi_t(\M).
\end{equation}
If all weights $w_i$ are equal, then \eqref{eq:MZ-curves} also holds without the factors $\frac{1}{3}$ and $3$. 
\end{proposition}
When $N$ is of order $t^d$, the curve $\gamma_w$ has length $\mathcal{O}(t^{d-1})$, matching the order achieved in \cite{EhlerGroechenigKarner} with non-simple curves for the sphere $\Sd$. 

The proof of Proposition \ref{prop:MZ-curves} is contained in Appendix \ref{app:MZ+LS}.

\subsection{Exact integration measures beyond points and curves}
The measure-theoretic formulation of Theorem \ref{thm:foundation} also allows us to construct exact integration measures of more general forms, including weighted combinations of locally supported probability measures and normalized restrictions of the Riemannian measure to unions of small balls. 

The following corollary serves primarily to illustrate the versatility of Theorem \ref{thm:foundation}. 
\begin{corollary}\label{cor:1}
Let $N$ be a positive integer and $R_1,\ldots,R_N$ an equal-area partition of $\M$ as in Theorem \ref{thm:partition}. If $\sigma_i$ is a probability measure supported on $R_i$, for $i=1,\dots,N$, then there are weights $a_i\in[\frac{1}{2N},\frac{3}{2N}]$ such that $\sigma=\sum_{i=1}^N a_i\sigma_i$ exactly integrates $\Pi_t(\M)$ for $t= \frac{1}{72A_{8b}}N^{1/d}$. 
\end{corollary}
\begin{proof}
For $L=1$ and $h=1/(2N)$, we have $\frac{L}{hN}=2\in[\tfrac{3}{2},\tfrac{5}{2}]$, and we define 
\begin{equation*}
\nu_{i,w}:= h w_i \sigma_i=\frac{w_i}{2N}\sigma_i,\qquad w\in[1,3]^N.
\end{equation*}
The oscillation estimate in Proposition \ref{lem:manifold-oscillation} yields, with $R_i\subseteq V_i$, for $P\in\Pi_t(\M)$ such that $\|P\|_{L^1(\M)}=1$, 
\begin{equation*}
 \frac{12}{N}\sum_{i=1}^N\osc_{V_i}(P)
  \leq 48 A_{8b} t N^{-1/d} \|P\|_{L^1(\mathcal M)} \leq \frac{48}{72} = \frac{2}{3}<1.
\end{equation*}
Theorem \ref{thm:foundation} implies that there are weights $\hat{w}\in [1,3]^N$ such that $\sigma:=\sum_{i=1}^N \frac{\hat{w}_i}{2N}\sigma_i$ integrates $\Pi_t(\M)$ exactly. We can therefore choose $a_i=\frac{\hat{w}_i}{2N}$.
\end{proof}
A more substantial consequence is that the method yields families of Logvinenko-Sereda sets as studied in \cite{MarzoCerda,Ortega,Wang}.

\begin{corollary}\label{cor:LS}
Let $N$ be a positive integer, $q_N=aN^{-1/d}$, and let $R_1,\dots,R_N$ be a partition with points $x_1,\ldots,x_N$ as in Theorem \ref{thm:partition}. 
\begin{itemize}
\item[(i)] Then, for every 
\begin{equation}\label{eq:L rho}
0<L\le
\frac{2}{3}N\min_{i=1,\dots,N}\mu(B_{q_N}(x_i)),
\end{equation}
there are radii \(0<r_i\le q_N\) such that $S_N:=\bigcup_{i=1}^{N}B_{r_i}(x_i)$
satisfies $\mu(S_N)=L$ and 
\begin{equation}\label{eq:id for sets}
\frac1 L\int_{S_N} P(x)\dd\mu(x)
=
\int_{\mathcal M}P(x)\dd\mu(x)
\qquad \forall\,P\in\Pi_t(\mathcal M),\quad t \leq  \frac{N^{1/d}}{72 A_{8b}},
\end{equation}
with $\frac{L}{2N}\leq \mu(B_{r_i}(x_i))
    \leq \frac{3L}{2N}$, for $i=1,\ldots,N$.

\item[(ii)] For $\varepsilon\in(0,1)$, there is a constant $c_{\M,\varepsilon}$ such that for each fixed sufficiently small $L>0$, there exists $N_0$, so that  the resulting sets $(S_N)_{N\geq N_0}$ in (i) satisfy, for all positive $t\leq c_{\M,\varepsilon} N^{1/d}$ and every $1\leq p\leq \infty$,
\begin{equation}\label{eq:MZ-sets}
     \tfrac{1}{2}(1-\varepsilon)\|P\|_{L^p(\M)}
    \le
    \|P\|_{L^p(\M,L^{-1}1_{S_N}\mu)}
    \le
     \tfrac{3}{2}(1+\varepsilon)\|P\|_{L^p(\M)},\quad \forall\, P\in\Pi_t(\M).
\end{equation}
\end{itemize}
\end{corollary}
A family of sets $(S_N)_{N\geq N_0}$ satisfying the Marcinkiewicz-Zygmund inequalities \eqref{eq:MZ-sets}, with constants uniform in $t$ but possibly depending on $p$, is called a Logvinenko-Sereda family in \cite{MarzoCerda,Ortega,Wang}. The sets constructed here additionally satisfy exact integration identities and the constants are independent of $p$. 

\begin{remark}\label{rem:density}
The uniform volume bounds for small geodesic balls give
constants $c_0,C_0>0$ such that, for all sufficiently large $N\geq N_0$,
\[
   \frac{c_0}{N}\leq  \mu(B_{aN^{-1/d}}(x)),
    \qquad
    \mu(B_{(a+b)N^{-1/d}}(x))\leq\frac{C_0}{N},
    \qquad x\in\M.
\]
Thus, for every fixed $0<L\leq 2c_0/3$, 
\[
    L\leq\frac{2c_0}{3}
    \leq\frac23N\min_{1\leq i\leq N}\mu(B_{q_N}(x_i)),
\]
so the assumed upper bound on $L$ is satisfied.

Fix such an $L$. For each $x\in\M$, choose $i$ with $x\in R_i$.
Since $R_i\subseteq B_{bN^{-1/d}}(x_i)$ and
$r_i\leq aN^{-1/d}$, we have $B_{r_i}(x_i)\subseteq B_{(a+b)N^{-1/d}}(x)$. 
Consequently,
\begin{equation}\label{eq:low and so}
    \frac{\mu\bigl(S_N\cap B_{(a+b)N^{-1/d}}(x)\bigr)}
         {\mu(B_{(a+b)N^{-1/d}}(x))}
    \geq \frac{L/(2N)}{C_0/N}
    =\frac{L}{2C_0}>0.
\end{equation}
Let \(T_N:=N^{1/d}/(72A_{8b})\). Since $(a+b)N^{-1/d}=\frac{a+b}{72A_{8b}T_N}$, 
the relative density bound \eqref{eq:low and so} is at spectral scale \(T_N^{-1}\). It yields the \(L^2(\M)\)-norm equivalences of \cite[Theorem 1.1]{Ortega}, and the \(L^p(\M)\)-norm equivalences of \cite[Theorem 1]{Wang} for \(1<p<\infty\). The proof in Appendix \ref{app:MZ+LS} covers all \(1\le p\le\infty\), with constants independent of \(p\).
\end{remark}

\appendix

\section{Proofs for the torus}\label{app:torus}
The $L^1$-Bernstein inequality on the torus $\T^d$ is classical \cite[Equation~(3)]{TotikVarga}.
\begin{lemma}[$L^1(\T^d)$-Bernstein inequality]\label{lemma:Bernstein T}
For every positive integer $t$, we have
\begin{equation*}
\|\partial_j P \|_{L^1(\T^d)} \leq 2\pi t \|P\|_{L^1(\T^d)},\qquad \forall\,P\in\mathcal{T}_t(\T^d). 
\end{equation*}
\end{lemma}

Let $R_1,\ldots,R_N$ be the uniform partition of
$[-\frac12,\frac12]^d$ into axis-parallel cubes of sidelength $h$ with $h^{-1}$ an even integer and $N=h^{-d}$.

\begin{proof}[Proof of Proposition \ref{lem:torus-cubical-oscillation}]
Let $u,v\in R_i$ be points at which $P$ attains its maximum and minimum on $R_i$. Since $R_i$ is convex, the segment
$z(s)=v+s(u-v)$, $0\leq s\leq1$, lies in $R_i$. Since $R_i$ is an axis-parallel cube of sidelength $h$,
we have $|u_j-v_j|\leq h$ for every $j$. Consequently,
\begin{align*}
 \osc_{R_i}(P) = P(u)-P(v)=\int_0^1\frac{d}{ds}P(z(s))\,ds   &=\sum_{j=1}^d(u_j-v_j) \int_0^1\partial_jP(z(s))\,ds\\
  &\leq
    \sum_{j=1}^d h
    \int_0^1|\partial_jP(z(s))|\,ds\\
  &\leq h\sum_{j=1}^d\sup_{x\in R_i}|\partial_jP(x)|.
\end{align*}
Since the unit sphere in the finite-dimensional space $\mathcal T_t(\T^d)$ is compact, the maximum 
\begin{equation*}
C:=\max_{0\neq Q\in \mathcal T_t(\T^d)} \frac{\frac1N\sum_{i=1}^N\operatorname{osc}_{R_i}(Q)}{\|Q\|_{L^1(\T^d)}}
\end{equation*}
exists, and let $Q\in\mathcal T_t(\T^d)$ be a maximizer with $\|Q\|_{L^1(\T^d)}=1$. 

For every \(x,y\in R_i\), the triangle inequality gives
\[
|\partial_jQ(y)|
\le |\partial_jQ(x)|+|\partial_jQ(y)-\partial_jQ(x)|
\le |\partial_jQ(x)|+\operatorname{osc}_{R_i}(\partial_jQ).
\]
We take the supremum over \(y\) and then average over \(x\in R_i\) to derive
\begin{equation*}
\sup_{R_i}|\partial_jQ|
\le N\int_{R_i}|\partial_jQ(x)|\,\dd x+\osc_{R_i}(\partial_jQ).
\end{equation*}

Since $\partial_jQ\in\mathcal T_t(\T^d)$, the definition of $C$ yields
\begin{align*}
\frac1N\sum_{i=1}^N\sup_{R_i}|\partial_jQ|
&\le
\sum_{i=1}^N\int_{R_i}|\partial_jQ(x)|\,dx
+\frac1N\sum_{i=1}^N\operatorname{osc}_{R_i}(\partial_jQ)\\
&\le
\|\partial_jQ\|_{L^1(\mathbb T^d)}
+C\|\partial_jQ\|_{L^1(\mathbb T^d)}\\
&=(1+C)\|\partial_jQ\|_{L^1(\mathbb T^d)}.
\end{align*}

The segment estimate at the start of the proof applies to $Q$
as well as to $P$:
\[
    \osc_{R_i}(Q)\le h\sum_{j=1}^d\sup_{R_i}|\partial_jQ|.
\]
Since $Q$ attains the maximum defining $C$ and
$\|Q\|_{L^1(\T^d)}=1$, summing gives
\[
    C=\frac1N\sum_{i=1}^N\osc_{R_i}(Q)
    \le\frac hN\sum_{j=1}^d\sum_{i=1}^N
    \sup_{R_i}|\partial_jQ|.
\]
Lemma~\ref{lemma:Bernstein T} yields
\[
    C
    \le h(1+C)\sum_{j=1}^d\|\partial_jQ\|_{L^1(\T^d)}
    \le 2\pi dht(1+C).
\]

Thus $C\le 2\pi dht/(1-2\pi dht)$ and $ht\leq\frac{1}{4\pi d}$ leads to $C\le4\pi dht$,
which proves Proposition \ref{lem:torus-cubical-oscillation}.
\end{proof}

\begin{proof}[Proof of Lemma \ref{lemma:parameters need balancing}]
The assumption on $L$ yields $L\ge 2(64\pi dt)^{d-1}$, so that 
\begin{equation*}
h^{-1}= 2\Big\lceil \frac{(L/2)^{1/(d-1)}}{2} \Big\rceil \geq 2\lceil 32\pi dt \rceil \geq 64 \pi d t.
\end{equation*}
Rearranging leads to $ht \leq 1/(64\pi d)$.

Set $r:=(L/2)^{1/(d-1)}$. Then $L=2r^{d-1}$ and $h^{-1}=2 \lceil\frac r2\rceil$ with $r\leq h^{-1}<r+2$. 
Consequently,
\[
  0<hr\leq1,
  \qquad
  0\leq1-hr
  =h\bigl(h^{-1}-r\bigr)\leq 2h.
\]
Since $N=h^{-d}$, we have $hN=h^{1-d}$.
The upper bound therefore follows from
\[
  \frac{L}{hN}
  =\frac{2r^{d-1}}{h^{1-d}}
  =2(hr)^{d-1}\leq2.
\]
For the lower bound, use $hr\leq1$ and the estimate
$th\leq 1/(64\pi d)$, so that we derive
\begin{align*}
  1-(hr)^{d-1}
  =(1-hr)\sum_{j=0}^{d-2}(hr)^j
&  \leq2h(d-1)\\
&   \leq2hd \leq\frac1{32\pi}
   <\frac14.
\end{align*}
Consequently, 
\[
  \frac{L}{hN}
  =2(hr)^{d-1}
  \geq2 (1-\tfrac14)
  =\frac32.
  \qedhere
\]
\end{proof}

\section{Embedded minimum spanning trees}\label{app:embedded}
We verify that the tree $T$ in Lemma \ref{lemma:embedded manifold} is embedded in three steps.

(i) First, no vertex can lie in the interior of an edge of $T$. Indeed, suppose that $x_k$
lies in the interior of the geodesic edge $e_{ij}$. Then
\[
    d_{\M}(x_i,x_j)
    =
    d_{\M}(x_i,x_k)+d_{\M}(x_k,x_j),
\]
and both terms on the right are strictly smaller than
$d_{\M}(x_i,x_j)$. Removing the edge $e_{ij}$ separates $T$ into two
components. According to which component contains $x_k$, one can
replace $e_{ij}$ by either $e_{ik}$ or $e_{kj}$ and obtain a spanning tree of
strictly smaller total length, contradicting the minimality of $T$.

(ii) Suppose distinct edges \(e_{ij}\) and \(e_{ik}\) share the vertex \(x_i\), and assume that they intersect at another point \(p\neq x_i\). Let \(g_j:[0,L_j]\to M\) and \(g_k:[0,L_k]\to M\) be the unit-speed minimizing geodesics from \(x_i\) to \(x_j\) and \(x_k\), respectively. Set \(s=d_M(x_i,p)\). Because both edges are minimizing,
\[
g_j(s)=g_k(s)=p,
\qquad 0<s<\min(L_j,L_k).
\]
Moreover, $s<h/2<\inj(\M)$. By uniqueness of the minimizing geodesic from \(x_i\) to \(p\),
\[
g_j|_{[0,s]}=g_k|_{[0,s]}.
\]
In particular, their velocities at \(p\) agree. Uniqueness for the geodesic initial-value problem therefore gives
\[
g_j(t)=g_k(t),
\qquad 0\leq t\leq\min(L_j,L_k).
\]
Since \(x_j\neq x_k\), we have \(L_j\neq L_k\). Consequently, the endpoint of the shorter edge lies in the interior of the longer edge, and case (i) applies.

(iii) Next, suppose that two nonadjacent edges $e_{i,j}$ and $e_{k,l}$ intersect in $p$. By (i), the point $p$ lies in the interior of both edges. If the tangent lines of the two edges at \(p\) coincide, uniqueness for the geodesic initial-value problem implies that the two geodesic segments are restrictions of the same geodesic. Since \(p\) is interior to both segments and their endpoints are distinct, at least one endpoint lies in the interior of the other segment, contradicting (i). Hence the two edges have distinct tangent lines at \(p\).
 
Starting with a tree, removing one edge after another increases the number of connected components by exactly one. 
Thus, removing $e_{i,j}$ and $e_{k,l}$ from $T$ produces three connected components. Neither the vertices $x_i$ and $x_j$ nor the pair $x_k$ and $x_l$ can be in the same connected component. 
After relabeling the vertices $x_i\leftrightarrow x_j$ and $x_k\leftrightarrow x_l$  if necessary, we may assume that
$x_i$ and $x_k$ belong to the same component, while $x_j$ and $x_l$ belong
to the other two components. Hence replacing $e_{i,j}$ and $e_{k,l}$ by
$e_{i,l}$ and $e_{j,k}$ again produces a spanning tree of $\widetilde{G}$.

Since the intersecting geodesics have distinct tangent lines at \(p\), both concatenations through \(p\) have a corner and are therefore not length-minimizing. Therefore, the triangle
inequality is strict at $p$, and we obtain
\begin{align*}
d_{\M}(x_i,x_l)
    &<
    d_{\M}(x_i,p)+d_{\M}(p,x_l),\\
    d_{\M}(x_j,x_k)
    &<
    d_{\M}(x_j,p)+d_{\M}(p,x_k).
\end{align*}
Since the original edges are minimizing geodesic segments and \(p\) lies on both, we have
\begin{align*}
d_\M(x_i,p)+d_\M(p,x_j)&=d_\M(x_i,x_j),\\
d_\M(x_k,p)+d_\M(p,x_\ell)&=d_\M(x_k,x_\ell).
\end{align*}
Adding the preceding strict inequalities therefore yields
\[
d_M(x_i,x_\ell)+d_M(x_j,x_k)
<
d_M(x_i,x_j)+d_M(x_k,x_\ell),
\]
contradicting the minimality of \(T\) in \(\widetilde G\). Thus the interiors of
distinct edges are disjoint.

\section{Initial curves in dimension two}\label{app:boundary}
We complete the proof of the case $d=2$ in
Proposition \ref{prop:simple initial curve}.
Recall that the tree $T$ is rooted at $x_1$.
For sufficiently small $\delta>0$, let
\[
D_{i,\delta}:=B_\delta(x_i),
\qquad i=1,\ldots,N,
\]
be the closed geodesic discs centered at the vertices of $T$.
Choose $\delta<\operatorname{inj}(M)$ so that these discs are
pairwise disjoint and each meets only the edges incident to
its centre.

Join the endpoint discs of each edge by a thin embedded
rectangular strip following that edge, as in the usual
construction of a neighborhood of an embedded graph
by vertex discs and edge bands; see
\cite[Sections~1.1.4--1.1.5]{EllisMonaghanMoffatt2013}.
Choose the strips pairwise disjoint, with their short ends
lying in the interiors of the corresponding endpoint discs,
and avoiding all other vertex discs.
Rounding the junctions between the discs and strips gives
a compact neighborhood $D_\delta$ of $T$ with smooth boundary.
Since a finite tree collapses to a point, $D_\delta$ is a
disc, see \cite[Corollary~3.27]{RourkeSanderson}.

Choose one of the two directions around the boundary $\partial D_{\delta}$ of
one sufficiently thin neighborhood $D_\delta$ of $T$.
Let $\Gamma:\mathbb{T}\to T$ be the corresponding tree
traversal, starting at $x_1$ and traversing every edge once
in each direction, in the order in which the two sides
of the edge strips of $\partial D_{\delta}$ are encountered, see \cite[Section~1.1, p.~249]{LeGall}. 
Put
\[
  \Lambda:=2\sum_{e \text{ edge of }T}\ell(e),
\]
and parametrize $\Gamma$ at constant speed $\Lambda$ on
$[0,1]$, with endpoints identified. Unlike the parametrization
used in the case $d\geq3$, this map need only be piecewise smooth.

Assign both traversals of each edge to its child, except that
the initial half, measured by arclength, of the first outgoing
traversal is assigned to the root. Let $I_1,\ldots,I_N$ be the
resulting finite unions of closed parameter intervals, which
partition $[0,1]$ up to shared endpoints. Every $I_i$ contains
a nondegenerate interval, and
\[
  \Gamma(I_i)\subseteq B_{h/4}(x_i),
  \qquad
  a_i:=\Lambda |I_i|\leq \frac h2.
\]
Indeed, each nonroot vertex is assigned at most the two
traversals of its parent edge, while the root is assigned
half of one edge traversal.

The neighborhoods can be chosen so that the two long sides
of each edge strip converge in $\mathcal C^1$, on compact
subarcs of the open edge, to that edge, and their lengths
converge to its length. The boundary arcs joining these sides
near the vertices, including the rounded leaf tips, have total
length tending to zero. Consequently,
\[
\ell(\partial D_\delta)\longrightarrow\Lambda,\qquad \text{when }\delta\to 0.
\]

Parametrize $\partial D_\delta$ by
$g_\delta:\mathbb{T}\to M$ at constant speed $\ell(\partial D_\delta)$,
so that the successive strip sides are traversed in the
order prescribed by $\Gamma$. 
Start at the root end of the strip side corresponding
to the first outgoing traversal of $\Gamma$.

The corresponding cumulative arclengths of the
successive edge sides converge to those of $\Gamma$.
Thus $g_\delta(s)\to\Gamma(s)$ for every $s\in\mathbb{T}$
such that $\Gamma(s)\notin\{x_1,\ldots,x_N\}$.
Only finitely many parameter values are excluded, since
$\Gamma$ traverses each edge exactly twice at constant speed.

The speeds $\ell(\partial D_\delta)$ are uniformly bounded, so the maps
$g_\delta$ are equicontinuous.
Together with the continuity of $\Gamma$, this implies
uniform convergence on all of $\mathbb{T}$:
\[
\sup_{s\in\mathbb{T}}
d_M\bigl(g_\delta(s),\Gamma(s)\bigr)\longrightarrow0,\qquad \text{when }\delta\to 0.
\]
Each $g_\delta$ is a smooth embedding with positive speed and
bounds the disc $D_\delta$.

Fix $0<\varepsilon<h/8$, and choose $\delta$ so small that
\[
  \sup_{s\in\T} d_\M(g_\delta(s),\Gamma(s))<\varepsilon,
  \qquad |\ell(\partial D_\delta)-\Lambda|<\varepsilon.
\]
Set $g=g_\delta$ and assign $g|_{I_i}$ to $x_i$.
This assigns the entire boundary, including all vertex
connections. Constant speed gives, for every $i$,
\[
  g(I_i)\subseteq B_{h/4+\varepsilon}(x_i),
  \qquad
  \int_{I_i}\|\dot g(s)\|\,ds
  =\ell(\partial D_\delta) |I_i|
  \leq a_i+\varepsilon|I_i|
  \leq\frac h2+\varepsilon.
\]
The interpolation and assignment arguments proceed as in the case \(d\geq 3\).

\section{Proof of Lemma \ref{lem:parametric-smoothing}}\label{app:corner}
Write $S=\mathbb R/m\mathbb Z$. We round the corners in normal
coordinates, using the same smoothing scale for every $w\in W$.

For each vertex $p_j(w)$, denote the incoming and outgoing velocities of
$\gamma_w$ by $v_{j,-}(w),v_{j,+}(w)$, respectively.
They are nonzero and depend continuously on $w$: indeed, the endpoints
of every edge vary continuously, and its length remains strictly below
$\inj(\M)$. In normal coordinates centered at $p_j(w)$, the two
incident geodesic segments are therefore represented, near $s=j$, by
\[
 q_{j,w}(s)=
 \begin{cases}
   (s-j)v_{j,-}(w),&s\le j,\\
   (s-j)v_{j,+}(w),&s\ge j.
 \end{cases}
\]

Choose a nonnegative, even mollifier
$\eta\in C_c^\infty((-1,1))$ with $_{-1}^1\int\eta(s)\dd s=1$ and put
$\eta_\delta(s)=\delta^{-1}\eta(s/\delta)$. For sufficiently small
$0<\delta<1/2$, define
\[
 q_{j,w,\delta}:=q_{j,w}*\eta_\delta.
\]
Since $q_{j,w}$ is linear on either side of $j$ and $\eta$ is even,
\[
 q_{j,w,\delta}(s)=q_{j,w}(s)
 \qquad\text{if }|s-j|\ge\delta.
\]
Replacing $\gamma_w$ on each $[j-\delta,j+\delta]$ by
$\exp_{p_j(w)}q_{j,w,\delta}$ thus defines a smooth closed curve
$\widetilde\gamma_{w,\delta}$ which agrees with $\gamma_w$ outside
these intervals. The replacements fit smoothly to all orders at their
endpoints. Because the exponential map is smooth and the incident
velocities vary continuously with $w$, the maps
\[
 W\ni w\longmapsto\widetilde\gamma_{w,\delta}\in \mathcal C^1(S,\M)
\]
are continuous for every fixed sufficiently small $\delta$.

(1) We next verify that a single sufficiently small $\delta$ makes every
$\widetilde\gamma_{w,\delta}$ embedded. Since the geodesic cycle $\gamma_w$ is simple,
the incoming and outgoing velocity vectors at a vertex cannot point
in opposite directions. Thus, 
\[
 u_{j,\pm}(w):=\frac{v_{j,\pm}(w)}{|v_{j,\pm}(w)|},
 \qquad z_j(w):=u_{j,-}(w)+u_{j,+}(w).
\]
yield $\langle u_{j,+}(w),u_{,j,-}(w)\rangle>-1$ and we get for either sign
\begin{align*}
\langle z_j(w),v_{j,\pm}(w)\rangle &= |v_{j,\pm}(w)| \langle u_{j,+}(w)+u_{j,-}(w),u_{j,\pm}(w)\rangle \\
& = |v_{j,\pm}(w)|(1+\langle u_{j,\mp}(w),u_{j,\pm}(w)\rangle)>0.
\end{align*}
The left-hand sides depend continuously on $w$, so their minimum over
all $j$ and $w\in W$ is a positive constant $c_0$.
Moreover, differentiation of the convolution gives
\[
 q'_{j,w,\delta}(s)
 =(1-\theta_{j,\delta}(s))v_{j,-}(w)
    +\theta_{j,\delta}(s)v_{j,+}(w),
 \qquad 0\le\theta_{j,\delta}(s)\le1.
\]
Consequently, we also have 
\[
 \big\langle z_j(w),q'_{j,w,\delta}(s)\big\rangle
 \ge c_0>0.
\]
The same inequality holds on the unchanged portions of the adjacent
geodesic edges in the normal chart. Thus the $z_j(w)$-coordinate is
strictly increasing through the entire rounded corner. In particular,
each rounded corner is injective and has nonvanishing tangent.

Next, we verify that the injectivity is uniform over $W$. Since $W$ is compact, we may 
choose a small parameter radius $\rho>0$, independent of $j,w$, such
that the original arcs on $[j-2\rho,j+2\rho]$ lie in the above normal
charts, and these parameter intervals are disjoint. For pairs of
parameters $s,s'\in S$ with distance $d_S(s,s')\geq \rho$, compactness
and the simplicity of every $\gamma_w$ give
\[
 r_*:=\min_{\substack{w\in W,\ s,s'\in S\\
                  d_S(s,s')\ge\rho}}
          d_{\M}(\gamma_w(s),\gamma_w(s'))>0.
\]
The incident edge lengths and the differentials of the exponential
maps are uniformly bounded on the chosen compact normal-coordinate
neighborhoods. Hence there exist $B>0$ and $\delta_0>0$, independent of
$w$, such that both $\gamma_w$ and
$\widetilde\gamma_{w,\delta}$ have speed at most $B$ for
$0<\delta\le\delta_0$. Comparing their positions with the closer
endpoint of each modified interval, where the curves coincide, yields
\begin{equation}\label{eq:param-smooth-C0}
 \sup_{w\in W,\ s\in S}
 d_{\M}(\widetilde\gamma_{w,\delta}(s),\gamma_w(s))
 \le2B\delta.
\end{equation}
Thus, when $4B\delta<r_*$, the modified curves cannot identify two
parameters at circular distance at least $\rho$. If two distinct
parameters are closer than $\rho$, then either they lie on one
unchanged minimizing-geodesic edge, or the short parameter arc between
them meets a rounded corner. In the latter case both parameters lie
in a common interval $[j-2\rho,j+2\rho]$, where the strictly
increasing $z_j(w)$-coordinate prevents coincidence. Therefore
$\widetilde\gamma_{w,\delta}$ is injective for all $w\in W$.
Its tangent is nonzero both at the rounded corners and on the
unchanged geodesic pieces; hence it is a smooth embedding.

If, in addition, $2B\delta<\inj(\M)$, then 
$\widetilde\gamma_{w,\delta}$ is homotopic to $\gamma_w$
by \cite[Problem~6-16(b)]{Lee}.

(2) Finally, the union of the modified parameter intervals has length
$2m\delta$. Outside this set the curves coincide, while the arclengths
of their two restrictions to it are each at most $2mB\delta$.
Therefore, for every $f\in \mathcal C(\M)$,
\begin{equation}\label{eq:param-smooth-integrals}
 \left|\int_{\tilde\gamma_{w,\delta}}f
              -\int_{\gamma_w}f\right|
 \le4mB\delta\,\|f\|_\infty,
 \qquad w\in W.
\end{equation}
Choose $\delta>0$ small enough for the embedding and homotopy
conditions above and so that
\[
 2B\delta\le\varepsilon,
 \qquad 4mB\delta\le\varepsilon.
\]
Then $\gamma_{w,\varepsilon}:=\widetilde\gamma_{w,\delta}$
has all the required properties.

\section{Proofs for Section \ref{sec:further}}\label{app:MZ+LS}

We first verify results on Marcinkiewicz-Zygmund inequalities. 
\begin{proof}[Proof of Proposition \ref{prop:MZ-curves}]
Put $V_i:=B_h(x_i)$ with $h=8bN^{-1/d}$. 
Then $R_i\subseteq V_i$ and $\gamma_{i,w}\subseteq V_i$.
The inner and outer ball conditions in Theorem \ref{thm:partition} imply that
the separation and covering radii of $\{x_1,\ldots,x_N\}$
are comparable to $N^{-1/d}$. The local gradient estimate in
\cite[Lemma~6.7, Equation~(6.22)]{FilbirMhaskar} therefore yields
\begin{equation}\label{eq:FM}
    \Big(
        \frac1N\sum_{i=1}^N
        \big(\sup_{x\in V_i}|\nabla_{\M} P(x)|\big)^p
    \Big)^{1/p}
    \le C t\|P\|_{L^p(\M)}
\end{equation}
whenever $t \le c N^{1/d} $ with the standard interpretation when $p=\infty$.
The constants $C$ and $c$ depend only on $\M$ and the fixed
partition constants $a$ and $b$, but are independent of $p$.

For every $y\in V_i=B_h(x_i)$, a minimizing geodesic from
$x_i$ to $y$ stays in $V_i$ and has length at most $h$.
Integration along this geodesic gives
\[
    |P(y)-P(x_i)|
    \le h\sup_{x\in V_i}|\nabla_\M P(x)|.
\]
Hence, by the triangle inequality,
\[
    \operatorname{osc}_{V_i}(P)
    =\sup_{y,z\in V_i}|P(y)-P(z)|
    \le 2h\sup_{x\in V_i}|\nabla_\M P(x)|.
\]
Consequently, with \eqref{eq:FM} 
we obtain
\begin{equation}\label{eq:MZ-local-oscillation}
    \Big(
            \frac1N\sum_{i=1}^N
            \operatorname{osc}_{V_i}(P)^p
        \Big)^{1/p}
    \le C_1t N^{-1/d} \|P\|_{L^p(\M)},
    \qquad 1\le p\le\infty,
\end{equation}
with $C_1$ independent of $p$. This generalizes Proposition \ref{lem:manifold-oscillation} from $p=1$ to its full range. 

Fix $w\in[1,3]^N$. Let $\nu_{i,w}$ be the arclength measure
induced by the assigned portion $\gamma_{i,w}$, and put
\[
    \nu_w:=\sum_{i=1}^N\nu_{i,w},
    \qquad
    \sigma_{i,w}:=\frac{\nu_{i,w}}{hw_i},
    \qquad
    \eta_w:=\frac1N\sum_{i=1}^N\sigma_{i,w}.
\]
Since $\nu_{i,w}(\M)=hw_i$, each $\sigma_{i,w}$ is a
probability measure supported in $V_i$.

Let $P\in\Pi_t(\M)$ and first suppose that $1\le p<\infty$.
Define the normalized local norms
\[
    a_i:=
    \Big(N\int_{R_i}|P(x)|^p\,\dd\mu(x)\Big)^{1/p},
    \qquad
    b_i:=
    \Big(\int_{\M}|P(x)|^p\,\dd\sigma_{i,w}(x)\Big)^{1/p}.
\]
Both $N1_{R_i}\mu $ and $\sigma_{i,w}$ are probability measures
concentrated in $V_i$. Consequently,
\[
    \min_{x\in V_i}|P(x)|
    \le a_i,b_i
    \le \max_{x\in V_i}|P(x)|,
\]
and therefore
\[
    |a_i-b_i|
    \le \operatorname{osc}_{V_i}(|P|)
    \le \operatorname{osc}_{V_i}(P).
\]
Moreover, since the $R_i$ form an equal-area partition,
\[
    \Big(\frac1N\sum_{i=1}^N a_i^p\Big)^{1/p}
    =\|P\|_{L^p(\M)},
    \qquad
    \Big(\frac1N\sum_{i=1}^N b_i^p\Big)^{1/p}
    =\|P\|_{L^p(\M,\eta_w)}.
\]
The reverse triangle inequality for the weighted
$\ell^p$-norm thus gives
\begin{align*}
    \left|
        \|P\|_{L^p(\M,\eta_w)}-\|P\|_{L^p(\M)}
    \right|
    &=
    \Big|
        \Big(\frac1N\sum_{i=1}^N b_i^p\Big)^{1/p}
        -
        \Big(\frac1N\sum_{i=1}^N a_i^p\Big)^{1/p}
    \,\Big|\\
    &\le
    \Big(\frac1N\sum_{i=1}^N|b_i-a_i|^p\Big)^{1/p}\\
    &\le
    \Big(
        \frac1N\sum_{i=1}^N
        \operatorname{osc}_{V_i}(P)^p
    \Big)^{1/p}.
\end{align*}
The case $p=\infty$ is treated analogously.

Choose $c_{\M,\varepsilon}:=\min (c,\varepsilon/C_1)$, so that $t  \le c_{\M,\varepsilon} N^{1/d}$ leads to $C_1t N^{-1/d}\leq \varepsilon$.  The estimate \eqref{eq:MZ-local-oscillation}
therefore implies
\begin{equation}\label{eq:MZ-balanced-measure}
    (1-\varepsilon)\|P\|_{L^p(\M)}
    \le
    \|P\|_{L^p(\M,\eta_w)}
    \le
    (1+\varepsilon)\|P\|_{L^p(\M)}.
\end{equation}
If all weights $w_i$ are equal, then $\eta_w= \frac{1}{\ell(\gamma_w)}\nu_w$ and \eqref{eq:MZ-balanced-measure} is what we need. If not, then  $1\le w_i\le3$ and
$hN\le \ell(\gamma_w)\le3hN$ yield $\frac13\eta_w
    \le
    \frac{\nu_w}{\ell(\gamma_w)}
    \le
    3\eta_w$. Hence, we derive 
\[
    \frac{(1-\varepsilon)}{3^{1/p}}\|P\|_{L^p(\M)}
    \le
    \|P\|_{L^p(\M,\gamma_w)}
    \le
    (1+\varepsilon) 3^{1/p}\|P\|_{L^p(\M)},
\]
where $3^{1/\infty}:=1$. The bounds in \eqref{eq:MZ-curves} follow.
\end{proof}

We verify the results on Logvinenko-Sereda sets.
\begin{proof}[Proof of Corollary \ref{cor:LS}]
(i) For $h=L/(2N)$, we have $\frac{L}{hN}=2\in[\tfrac{3}{2},\tfrac{5}{2}]$. The assumption \eqref{eq:L rho} ensures 
\begin{equation*}
\mu(B_{q_N}(x_i))\geq 3 \frac{L}{2N}=3h\geq h w_i.
\end{equation*}
Hence, there are maps $r_i:[1,3]^N\to (0,q_N]$ such that $\mu\big(B_{r_i(w)}(x_i) \big)=h w_i$, for $i=1,\dots,N$, 
and we may define 
\begin{equation*}
\nu_{i,w}:= 1_{B_{r_i(w)}(x_i)}\,\mu,\qquad i=1,\dots,N,\qquad \nu_w:=\sum_{i=1}^N \nu_{i,w}.
\end{equation*}
The continuity requirement follows from 
\begin{equation*}
\big|\int_{\M} f(x) \dd\nu_{i,w}-\int_{\M} f(x) \dd\nu_{i,w'} \big|\leq \|f\|_{L^\infty(\M)}  h |w_i-w_i'|,\qquad f\in\mathcal{C}(\M). 
\end{equation*}

We have $B_{q_N}(x_i)\subseteq R_i$. For $T_N:=\frac{N^{1/d}}{72 A_{8b}}$ and the normalized polynomials in \(\Pi_{T_N}(\mathcal M)\), Proposition \ref{lem:manifold-oscillation} yields $\frac{12}{N}\sum_{i=1}^N\osc_{V_i}(P)<1$ as in the proof of Corollary \ref{cor:1}. Theorem \ref{thm:foundation} implies that there is $\hat{w}\in[1,3]^N$ such that $\frac{1}{L}\nu_{\hat{w}}=\frac{1}{L}\sum_{i=1}^N 1_{B_{r_i(\hat{w})}(x_i)}\,\mu$ exactly integrates $\Pi_{T_N}(\M)$. 
Thus, $S_N:=\bigcup_{i=1}^{N}B_{r_i(\hat{w})}(x_i)$ satisfies \eqref{eq:id for sets} and $P\equiv 1$ yields $\mu(S_N)=L$.

The volume bounds follow from
\[
    \mu(B_{r_i}(x_i))
    =h\widehat w_i=\frac{L}{2N}\widehat w_i,
    \qquad \widehat w_i\in[1,3].
\]

(ii) We choose $L$ as in Remark \ref{rem:density}. The uniform oscillation estimate
in the proof of Proposition \ref{prop:MZ-curves} gives, for every $1\le p\le\infty$,
\[
\Big(
\frac1N\sum_{i=1}^N
\bigl(\operatorname{osc}_{V_i}P\bigr)^p
\Big)^{1/p}
\le \varepsilon \|P\|_{L^p(\M)},
\qquad
P\in\Pi_{t}(\M),\quad t\leq c_{\M,\varepsilon}N^{1/d},
\]
where \(V_i=B_{8bN^{-1/d}}(x_i)\), and the left-hand side is interpreted
as a maximum when \(p=\infty\). 

Write \(B_i=B_{r_i}(x_i)\) and \(m_i=\mu(B_i)\), and define the
probability measures
\[
\eta_N=\frac1N\sum_{i=1}^N\frac{1_{B_i}\mu }{m_i},
\qquad
\omega_N=\frac{1_{S_N}\mu}{L}.
\]
We proceed as in the proof of Proposition \ref{prop:MZ-curves}. The relations 
$\frac{L}{2N}\le m_i\le\frac{3L}{2N}$ yield $\frac12\eta_N\le\omega_N\le\frac32\eta_N$, which leads to 
\[
2^{-1/p}(1-\varepsilon)\|P\|_{L^p(\M)}
\le
\left(\frac1L\int_{S_N}|P(x)|^p\,d\mu(x)\right)^{1/p}
\le
\left(\frac32\right)^{1/p}(1+\varepsilon)\|P\|_{L^p(\M)},
\]
for $1\le p<\infty$. In the case \(p=\infty\), the measures \(\eta_N\) and \(\omega_N\) have
the same null sets, so their \(L^\infty\)-norms coincide. 
\end{proof}

\section*{Acknowledgement of AI assistance}
ChatGPT (OpenAI; model: GPT-6 Astra) was used to assist with the development and discussion of selected proof steps within a proof strategy specified by the author, as well as with language editing. The author takes full responsibility for the content of this paper.

\section*{Acknowledgments}
The author thanks Karlheinz Gr\"ochenig for valuable discussions and helpful suggestions.

\bibliographystyle{abbrv}

\end{document}